\documentclass[1p]{elsarticle}
\pdfoutput=1

\usepackage{fullpage}
\usepackage{tikz}
\usepackage{wrapfig}
\usepackage{bmpsize}
\usepackage{amsmath,bm,bbm,amsthm, amssymb}
\usepackage{color}
\usepackage{subfig}
\usepackage{tabularx}
\usepackage{multirow} 
\usepackage{pgf,tikz}
\usepackage{mathrsfs}
\usetikzlibrary{arrows}
\usepackage[hyperindex,breaklinks]{hyperref}
\hypersetup{ colorlinks=true, linkcolor=blue, citecolor=blue,filecolor=blue, urlcolor=blue}
\usepackage{enumitem}
\usepackage{dsfont}

\definecolor{wiasblue}   {cmyk}{1.0, 0.60, 0, 0}
\definecolor{mlugreen}{RGB}{172,6,52}

\makeatletter
 \def\paragraph{\@startsection{paragraph}{4}%
 \z@\z@{-\fontdimen2\font}%
   {\normalfont\itshape}}\makeatother

\def\Z{\mathbb Z}

\def\E{\mathbb E}
\def\P{\mathbb P}

\def\R{\mathbb R}
\def\mc{\mathcal}
\def\ms{\mathsf}

\def\la{\lambda}

\def\s{\sigma}
\def\su{\subseteq}

\def\e{\varepsilon}
\def\t{\tau}

\def\b{\beta}
\def\de{\delta}
\def\et{\eta}

\def\es{\varnothing}
\def\one{\mathbbmss{1}}

\def\De{\Delta}
\def\G{\Gamma}

\def\ff{\infty}

\def\vp{\varphi}

\def\d{{\rm d}}

\def\AA{\mc A}
\def\BB{\mc B}
\def\CC{\mc C}

\def\GG{\mc G}

\def\NN{\mc N}
\def\PP{\mc P}
\def\XX{\mc X}

\def\TT{\mc T}
\def\UUU{\ms U}

\def\f{\frac}
\def\lan{\langle}
\def\ran{\rangle}

\def\r{\rho}

\def\im{\item}
\def\sm{\setminus}

\def\bep{\begin{proof}}
\def\enp{\end{proof}}
\def\bepr{\begin{proposition}}
\def\enpr{\end{proposition}}
\def\bec{\begin{corollary}}
\def\enc{\end{corollary}}
\def\bea{\begin{align}}
\newcommand\eea{\end{align}}
\def\beas{\begin{align*}}
\def\eeas{\end{align*}}
\def\bet{\begin{theorem}}
\def\ent{\end{theorem}}
\def\bee{\begin{example}}
\def\ene{\end{example}}

\def\bede{\begin{definition}}
\def\ende{\end{definition}}
\def\ber{\begin{remark}}
\def\enr{\end{remark}}
\def\beca{\begin{cases}}
\def\enca{\end{cases}}
\def\bel{\begin{lemma}}
\def\enl{\end{lemma}}
\def\been{\begin{enumerate}}
\def\enen{\end{enumerate}}

\def\beit{\begin{itemize}}
\def\enit{\end{itemize}}
\def\befr{\begin{frame}}
\def\enfr{\end{frame}}
\def\ti{\times}

\def\Var{\ms{Var}}
\def\Cov{\ms{Cov}}
\def\tf{\tfrac}
\def\ba{|\,}

\def\bi{\big}

\def\diam{\ms{diam}}

\def\wt{\widetilde}

\renewcommand\le{\leqslant}
\renewcommand\ge{\geqslant}

\def\dist{\ms{dist}}
\def\Rip{\ms{Rip}}

\def\becb{\begin{tcolorbox}[colback=Dandelion!20]}
\def\encb{\end{tcolorbox}}

\def\co{\colon}
\def\bn{\b^{M, b, d}_n}
\def\bnn{\b^{b, d}_n}

\def\PD{\ms{PD}}
\def\Eb{E_{\ms b}}
\def\Ed{E_{\ms d}}
\def\Eh{E_{\ms h}}
\def\Ebd{E_{\ms {bd}}}
\def\Var{\ms{Var}}
\def\kn{k_n}
\def\rh{\rho_{h, q, n}}
\def\dnz{\De_{i, n}}
\def\dex{\De_{z_i, n}'(h)}
\def\dexx{\De_{z_i, n}'}
\def\Qn{Q_n}
\def\TTP{T_{\ms{TP}}}
\def\TM{T_{\ms M}}

\theoremstyle{plain}
\newtheorem{theorem}{Theorem}[section]
\newtheorem{proposition}[theorem]{Proposition}
\newtheorem{corollary}[theorem]{Corollary}
\newtheorem{lemma}[theorem]{Lemma}

\theoremstyle{definition}
\newtheorem{definition}[theorem]{Definition}

\newtheorem{example}[theorem]{Example}

\theoremstyle{remark}
\newtheorem{remark}[theorem]{Remark}

\journal{\textcolor{white}{\; }}

\begin{document}

\begin{frontmatter}

	\title{Topology-based goodness-of-fit tests for sliced spatial data}

	\author[1]{Alessandra Cipriani}
\ead{A.Cipriani@tudelft.nl}
\address[1]{TU Delft (DIAM), Building 36, Mekelweg 4, 2628 CD, Delft, The Netherlands}
	\author[2,3,4]{Christian Hirsch}
\ead{hirsch@math.au.dk}
	\address[2]{Department of Mathematics,  Aarhus University,   Ny Munkegade 118, 8000, Aarhus C,  Denmark}
	\address[3]{Bernoulli Institute for Mathematics, Computer Science and Artificial Intelligence,  University of Groningen,  Nijenborgh~9, NL-9747 AG Groningen, The Netherlands}
	\address[4]{CogniGron (Groningen Cognitive Systems and Materials Center),  University of Groningen,  Nijenborgh 4, NL-9747 AG Groningen, The Netherlands}
	\author[1,5]{Martina Vittorietti}
\ead{M.Vittorietti@tudelft.nl}
\address[5]{ Department of Economics, Business and Statistics, University of Palermo, Viale delle Scienze, Building 13, Palermo, Italy}

\begin{abstract}
	In materials science and many other application domains, 3D information can often only be extrapolated by taking 2D slices. In topological data analysis, persistence vineyards have emerged as a powerful tool to take into account topological features stretching over several slices. In the present paper, we illustrate how persistence vineyards can be used to design rigorous statistical hypothesis tests for 3D microstructure models based on data from 2D slices. More precisely, by establishing the asymptotic normality of suitable longitudinal and cross-sectional summary statistics, we devise goodness-of-fit tests that become asymptotically exact in large sampling windows. We illustrate the testing methodology through a detailed simulation study and provide a prototypical example from materials science.
\end{abstract}


\end{frontmatter}


\section{Introduction}
\label{int_sec}

Topological data analysis (TDA) is an emerging branch within the domain of data science that holds the promise to unearth subtle properties of data by extracting shape-related characteristics. While TDA has its root in algebraic topology, which is often considered as one of the most theoretical fields of mathematical research, it is now applied in a variety of application domains such as astronomy, biology and materials science \cite{wasserman}.

While extracting refined topological information is already highly non-trivial, this becomes even more difficult when the data varies over time. To address this challenge, the notion of \emph{persistence vineyards} was developed, which makes it possible to track the evolution of topological features over time \cite{vine}. These time-varying features are also known as \emph{vines}.
While the initial motivation for vineyards arose in a problem of protein folding, this method has now been successfully used in a variety of further application contexts such as brain activity data gained via electroencephalography and functional magnetic resonance imaging \cite{fmri,brain}. In essence, vineyards allow to synthesize the information inherent in different time slices into a global summary.

In this manuscript, we argue that the concept of vineyards is not by any means restricted to time-varying data but extends to far more general situations involving sliced data. As a prototypical example, we highlight a data set from materials science where the 3D microstructure of an extra low carbon strip steel is measured via 2D slices. Developing statistical tests to determine whether a given stochastic-geometry model provides a good fit to the given sample is an important task since the properties and performance of metals are intimately linked to its microstructure. In fact, unravelling the relation between microstructure features and mechanical properties can lead to the design of new materials with desired properties \cite{stiff}.

The first step in the investigation is material characterization.
The most common way to observe microstructures is by 2D characterization techniques such as light microscopy, scanning electron microscopy, and electron backscatter diffraction. For opaque materials, serial slicing can be used for direct quantitative characterizations of 3D microstructures \cite{tewari}. Conventional serial slicing involves photographing (or digitally recording) a microstructural field-of-view, polishing the specimen to remove small thickness and, in the second metallographic plane, photographing the field-of-view exactly below the first one. This procedure (i.e., photograph-polish-etch-photograph) is then repeated to generate a stack of large number of aligned serial slices from which 3D microstructure information can be obtained. 

However, this process causes material disruption and is time-consuming especially if misalignment problems occur \cite{pirgazi}. Finding a stochastic model for microstructures is therefore critical in the study of the relation between microstructures features and mechanical properties. There are also important conceptual contributions towards the reconstruction of 3D structures from 2D slices such as \cite{amini}. However, in practice, it is difficult to assess to what extent the proposed density and transversality constraints are satisfied for a given dataset. Moreover, implementing the suggested method is not straightforward. Due to these difficulties, the current testing approaches for 3D data focus almost exclusively on test statistics computed from isolated 2D slices \cite{hahn,vittorietti}. Therefore, we are in urgent need for tests that can detect deviations of a model from the actual 3D topology without going through the process of a fully-fledged reconstruction.

Motivated from this example, we present a flexible testing framework for statistics derived from persistence vineyards. In its general form, this framework encompasses averages over certain scores associated with the individual vines. However, a large spacing between slices may cause uncertainties when identifying the vines. Therefore, we also propose a simpler class of cross-sectional test statistics that can be computed entirely from aggregations over the single slices. Our goodness-of-fit-tests are based on the asymptotic normality of these test statistics in large domains under suitable stabilization conditions. This asymptotic normality is the main conceptual contribution from our work. Our framework takes into account that when working with data from materials science, it is quite common to find samples in a highly imbalanced sampling window whose extension in the $x$-$y$-directions is far larger than in the $z$-direction.   We stress that the asymptotic normality of the test statistics in large windows offers two decisive advantages. First, computing confidence intervals only requires the mean and variance under the null model. Second, it suffices to compute these quantities for moderately large sampling windows; the values obtained in this manner can then be used for data coming from arbitrarily large sampling windows.

After having established our testing framework, we illustrate the method in a simulation study of 2D slices taken from a 3D Voronoi tessellation. This setup illustrates that the asymptotic normality is already clearly visible for moderately large sampling windows. Moreover, by considering tessellations induced by different classes of point processes as generators, we provide indications on the statistical power of the testing framework.
Finally, we apply the proposed method to the highlighted example of a 3D metallic microstructure. In this part, we also elucidate which practical obstacles need to be addressed when working with real datasets.

The rest of the manuscript is organized as follows. In Section \ref{mod_sec}, we define precisely our stochastic-geometry model and the considered test statistics. We also state our results on the asymptotic normality of the test statistics in large domains and explain the sufficient conditions.  Among all the possible models used for representing materials microstructure, Voronoi tessellations stand out since even this basic model  proved its power in approximating single-phase microstructures. The attractive mathematical properties and the availability of a wide range of sub-models make Voronoi tessellations the state of art for modelling microstructures \cite{madej}. Therefore, we verify in Section \ref{vor_sec} that the conditions of our main results are satisfied for the key example of 2D sections from a 3D Poisson-Voronoi tessellation. Next, in Section \ref{simulation_sec}, we illustrate the proposed testing methodology on several simulated data sets. In Section \ref{data_sec}, we analyze 2D sections taken from a sample of an extra low carbon strip steel. The proofs for the asymptotic normality of the test statistics are presented in the supplementary material. Finally, Section \ref{conc_sec} summarizes the findings and provides an outlook to future research.

\section{Methodology and main results}
\label{mod_sec}
In this section, we explain the general testing methodology and present the main results. First, in Section \ref{modd_sec}, we describe the model assumption of a point cloud varying between sections. Second, in Section \ref{pers_sec}, we present the $M$-bounded persistence diagram as the main tool to extract topology-related information from data. Next, Section \ref{stat_sec} describes a class of specific test statistics derived from the persistence diagram. Finally, Section \ref{norm_sec} contains the precise statement of the results on the asymptotic normality of these test statistics provided that a crucial set of stabilization conditions, which are described in detail in Section \ref{cond_sec}, is satisfied.

\subsection{Model }
\label{modd_sec}

We consider statistical testing problems for data describing some geometric structure in the $p$-dimensional Euclidean space, which is measured along $p'$-dimensional slices with $p' < p$. In the example from materials science presented in Section \ref{data_sec}, we consider a 3D structure measured via 2D slices, i.e., $p = 3$.. However, for clarity of exposition, we present the overall framework in general dimensions.

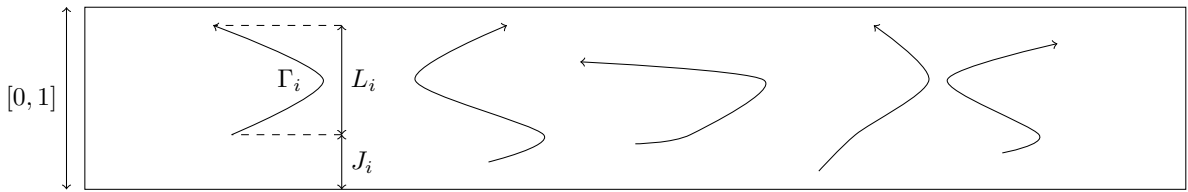
\begin{figure}[!htpb]
	\centering
	\resizebox{.99\linewidth}{!}{\begin{tikzpicture}[scale = 2.5]

\draw (-3, 0) rectangle (3,1);

\draw[->]  plot [smooth] coordinates {  (-2.2,.3) (-1.7,.6) (-2.3,.9)};
\draw[->]  plot [smooth] coordinates { (-0.8,.15) (-0.5,.3) (-1.2,.6) (-0.7,.9)};
\draw[->]  plot [smooth] coordinates { (0,.25) (0.3,.3) (0.7,.6) (-0.3,.7)};
\draw[->]  plot [smooth] coordinates { (1,.1) (1.2,.3) (1.6,.6) (1.3,.9)};
\draw[->]  plot [smooth] coordinates { (2,.2) (2.2,.3) (1.7,.6) (2.3,.8)};

\draw[<->] (-1.6, .3)--(-1.6,.9);
\draw[<->] (-1.6, .3)--(-1.6,.0);
\draw[dashed] (-1.6, .3)--(-2.2, .3);
\draw[dashed] (-1.6, .9)--(-2.3, .9);
\draw[<->] (-3.1, .0)--(-3.1,1);

\coordinate[label=0: {$L_i$}] (D) at (-1.6, .6);
\coordinate[label=0: {$J_i$}] (D) at (-1.6, .15);
\coordinate[label=0: {$\Gamma_i$}] (D) at (-2.0, .6);
\coordinate[label=-180: {$[0, 1]$}] (D) at (-3.1, .5);

\end{tikzpicture}}
	\caption{Point process of trajectories within the window $Q_n$, $p'=1$}
	\label{traj_fig}
\end{figure}

%
%
More precisely, in the null model, we let  $\XX_n = \{X_{i, n}(\cdot)\}_{i \ge 1}$ be a process of trajectories in a sampling window $Q_n := Q_n' \times [0, 1]$,  where $Q_n' := [-n/2, n/2]^{p'}$ is the $p'$-dimensional cube of side length $n$. See Figure \ref{traj_fig} for an illustration in 2D.  Hence, we always think of the last coordinate as the one along which slices are taken. In the data example discussed in Section \ref{data_sec}, the material is represented as a 3D tessellation, and the trajectories are formed by the centroids of the 2D cells in the slices.

On the modeling side, the $i$th trajectory is determined by an \emph{offset} $J_{i, n} \in [0, 1]$, a \emph{length} $L_{i, n} \in [0, 1]$ and  a \emph{shape} $\G_{i, n}$, where the latter is an element of $C([0, 1], Q_n)$, the space of continuous curves from $[0, 1]$ to $Q_n$. More precisely, for $\{(J_{i, n}, L_{i, n}, \G_{i, n})\}_{i \ge 1} \su [0, 1] \times [0, 1] \times C([0, 1], Q_n)$ we set 
$$X_{i, n}(h) := \big\{\big( \G_{i, n}((h - J_{i, n}) / L_{i, n}), h\big) \co J_{i, n} \le h \le  J_{i, n} + L_{i, n}\big\}.$$ 

We point out that when working with datasets stemming from a finite number of slices, some effort is needed for the trajectory reconstruction. We will discuss this issue in greater detail in Section \ref{data_sec}.

%
%
\subsection{Method: $M$-bounded persistence diagrams}
\label{pers_sec}
Having set up the model, we rely on the tools of \emph{persistent homology} in order to quantify the most striking topological aspects of the slices. For a superbly-written overview of this methodology intended for an audience with a background in statistics, we suggest \cite{wasserman} for further reading. However, for the presentation of the main results we do not need this abstract setting since we focus on the case $p = 3$ from now on. The restriction to the case of 2D slices of a 3D structure has the following reasons:
\been
\im The data example in Section \ref{data_sec}, which provides the motivation for the development of our goodness-of-fit tests, pertains precisely to the scenario of 2D slices of a 3D material. 
\im When working with 2D slices, then persistent homology captures two types of topological features, namely connected components and two-dimensional holes. As we will elaborate further below, both of these features admit a concrete topological description not relying on the abstract toolkit of simplicial homology.
\im Long-range correlations coming from percolation theory render it difficult to keep track of features of an arbitrary size. Therefore, we proceed as in \cite{svane} and restrict our attention to features whose size admits a deterministic upper bound. Again, while in general it is not at all immediate how to measure the size of features represented in simplicial homology, the 2D setting has already been worked out in detail in \cite{svane}.
\enen

Let $M > 0$ and $\vp \su \R^2$ be a finite set of points. In order to make the manuscript self-contained, we now briefly recall the concepts of $M$-bounded clusters and $M$-bounded holes (or loops), and refer the reader to \cite[Section 2]{svane} for the full technical details.

\begin{figure}[!h]
	\centering
\begin{tikzpicture}[scale=-1.0]
\draw[fill, black!29] (.01,0) circle (0.75);
\draw[fill, black!29] (-1.01,-0.8) circle (0.75);
\node at (-0.21,0.1) {$x$};
\node at (-2.21,-1.6) {\text{\bf a.}};
\draw[fill, black!29] (1.81,0.9) circle (0.75);
\draw[fill] (.01,0) circle (0.05);
\node at (-1.21,-0.7) {$z$};
\draw[fill] (-1.01,-0.8) circle (0.05);
\node at (1.61,1) {$y$};
\draw[fill] (1.81,0.9) circle (0.05);
\end{tikzpicture}
\qquad
\begin{tikzpicture}[scale=-1.0]
\draw[fill, black!29] (-0.72,0.7) circle (0.75);
\draw[fill, black!29] (0.72,-0.8) circle (0.75);
\draw[fill, black!29] (.01,1) circle (0.75);
\draw[fill, black!29] (1.01,0) circle (0.75);
\draw[fill, black!29] (-0.7101,-0.71) circle (0.75);
\draw[fill, black!29] (2.01,0.3) circle (0.75);
\draw[fill, black!29] (-1.01,0) circle (0.75);
\draw[fill] (-0.72,0.7) circle (0.05);
\draw[fill] (0.72,-0.8) circle (0.05);
\draw[fill] (.01,1) circle (0.05);
\draw[fill] (1.01,0) circle (0.05);
\draw[fill] (-0.7101,-0.71) circle (0.05);
\draw[fill] (2.01,0.3) circle (0.05);
\draw[fill] (-1.01,0) circle (0.05);
\draw (-0.72,0.7) --(.01,1)-- (1.01,0) -- (0.72,-0.8) -- (-0.7101,-0.71)-- (-1.01, 0) -- (-0.72,0.7);
\draw[fill] (.01,0) circle (0.05);
\node at (-0.531,0.0) {$p(H)$};
\node at (-2.201,-1.6) {\text{\bf b.}};
\end{tikzpicture}
	\caption{a. M-bounded clusters, b. M-bounded holes}
	\label{mb_fig}
\end{figure}
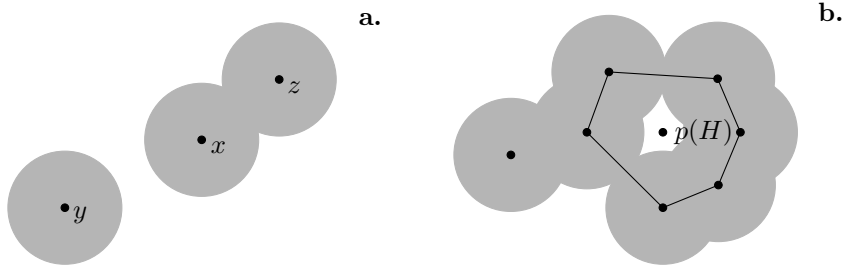

\subsubsection{$M$-bounded clusters}
 To describe the $M$-bounded clusters at a level $r > 0$, we consider the connected components of the union set $\bigcup_{x \in \vp}B_r(x)$, where $B_r(x)$ denotes the Euclidean disk of radius $r > 0$ centered at $x \in \vp$.  Note that at level $r = 0$ each $x \in \vp$ forms its own connected component. We say that all of these components are born at $r = 0$. Each of these components disappears at a unique level $r$, which is called the \emph{death time} of the component. This may occur for one of two reasons: i) at the level $r$ the spatial diameter of the component exceeds $M$ or ii) at the level $r$ two growing components merge into a single one. Here, we stress that in case ii) only one of the two components ``dies''. For instance, we may proceed as in \cite[Section 2.1]{svane} and decide this according to the lexicographic ordering of the two disks from the components that collide at level $r$. See Figure \ref{mb_fig} for an illustration (adapted from \cite{svane}).

The collection of death times $\{D_i\}_{i \ge 1}$ defines the \emph{$M$-bounded persistence diagram in degree 0}, which gives a summary of the $M$-bounded 0-features, i.e., of the merging pattern of $M$-clusters. To be consistent with the persistence diagram of $M$-bounded features to be introduced below, sometimes an additional first coordinate is added, which is identical to 0. We then write $\{(0, D_i)\}_{i \ge 1}$.

\subsubsection{$M$-bounded holes}
To describe the $M$-bounded holes, we proceed along similar lines with the difference that now we do not consider the connected components of $\bigcup_{x \in \vp}B_r(x)$ but rather the connected components of the complement $\R^2 \sm \bigcup_{x \in \vp}B_r(x)$ that are of size at most $M$. To any such component $H$ we may assign the \emph{center $p(H)$} to be the point which is covered last, i.e., which has the largest distance to $\vp$. If that point is covered at a level $r > 0$, then we refer to that level as the \emph{death time} of that component.  In contrast to the situation of the clusters, there are no $M$-bounded holes at level $r = 0$. A new $M$-bounded hole may arise at a level $r > 0$ if an existing (bounded or unbounded) component of $\R^2 \sm \bigcup_{x \in \vp}B_r(x)$ is split up into two. Next, we expound on how to measure the size of holes.  More precisely, $B_r(x_1), \dots, B_r(x_k)$ are the disks defining the boundary of a bounded connected component 
of $\R^2 \sm \bigcup_{x \in \vp}B_r(x)$, then we define the \emph{size} of that component as the metric diameter of the set $\{x_1, \dots, x_k\}$. An \emph{$M$-bounded hole} is a connected component of size at most $M$. Finally, we define the \emph{birth time} of a hole $H$ to be the smallest level $r > 0$ for which there exists an $M$-bounded hole $H'$ with $p(H') = p(H)$. 

Similarly to the setting of $M$-bounded clusters, the collection of $M$-bounded birth times and death times $\{(B_i, D_i)\}_{i \ge 1}$ is a fingerprint for the appearance and disappearance of $M$-bounded 1-features. Hence, $\{(B_i, D_i)\}_{i \ge 1}$ is called the \emph{$M$-bounded persistence diagram in degree 1}. In particular, both the birth time and the death time of $M$-bounded features admit a deterministic upper bound, and we may therefore assume that the persistence diagram is contained in $[0, \t]^2$ for some $\t > 0$.

\subsection{Test statistics}
\label{stat_sec}
The core of our testing methodology may be summarized in two key steps. First,  in Section \ref{pers_sec}, we  leveraged  persistence diagrams to extract topological information from each of the slices. Now, we describe how to combine the information from individual slices into more powerful global test statistics.

To that end, in each slice $h \in [0, 1]$, we compute the $M$-bounded persistence diagram induced by points $\{X_{j, n}(h)\}_{j \ge 1}$, tacitly discarding trajectories with $h \not \in [J_{j, n}, J_{j, n} + L_{j, n}]$. The persistence diagram in dimension $q$ tracks the appearance and disappearance of $q$-dimensional features. Thus, it is represented as a collection  $\{(B_i(h), D_i(h))\}_i = \{(B_{i, n}^q(h), D_{i, n}^q(h))\}_i$
of birth and death times in a fixed time interval $[0, \t]$.

%
%
Since the persistence diagrams in different slices are intimately related through the trajectories of the underlying point pattern, it is attractive to move beyond considering the respective persistence diagrams in isolation. As indicated in Section \ref{int_sec}, persistence vineyards or time-varying persistence diagrams are introduced in \cite{vine} as a tool to study the persistence of dynamic point clouds. Although the vines as described in the previous paragraph are similar in spirit to this classical concept, our setting deviates from the time-varying vines in a few technical characteristics, on which we now elaborate in further detail. First, \cite[Section 4]{vine} computes the simplicial complex based on the lower-star filtration of an embedding of an abstract set of vertices. Examples considered in the literature comprise distances resulting from the embedding a pair of points as single vertex, or the marking of 1-simplices through the negative correlation of two time series.  In contrast, the vertices in our vines are already embedded in Euclidean space so that we may rely on the classical \v Cech filtration. Even though with a sufficiently fine spacing between the slices, we are in principle able to track the vertices over the slices, we still observe \emph{switches} as in  \cite[Chapter VIII]{edHar}. That is, when the filtration values of certain simplices swap their ordering, then the simplices causing the births and deaths of features may change from one slice to another.

After stacking the persistence diagrams along the $h$-axis, we now elucidate how to define a \emph{vine} $(B_i(h), D_i(h)), h \in [0, 1]$ associated with a \emph{fixed feature $i$.} 
 Each point in the diagrams traces out a path called a \textit{vine}, i.e., a line among points of the persistence diagrams. We start with the case of $0$-features, i.e., $M$-bounded clusters. Here, the features are in one-to-one correspondence with the original trajectory process. That is, for $i\ge1$ and $h \in [J_i, J_i + L_i]$, we define $D_i(h)$ to be the death time of the point $X_{i,n}(h)$ on $\{X_{j, n}(h)\}_{j \ge 1}$. 

The case of 1-features is more involved since this one-to-one correspondence breaks down. In this case, we associate the $i$th feature with the simplex leading to the birth of that feature. Note that whether or not a certain collection of centers defines a feature may change from one slice to another. We stress that this convention only leads to a certain approximation of the actual vines. Moreover, it may often be very natural to consider a hole as being associated with the same feature even if a different simplex causes the birth of this hole. Finding useful heuristics for carrying out these more subtle identifications is an instance of the problem of \emph{cycle registration} and is an active research field \cite{reani}.

%
 %

%
%
To test whether data gathered from slices comply with a null hypothesis on the distribution on the underlying full-dimensional geometric structure, we propose the following general class of \emph{longitudinal test statistics}. In analogy to the vines from \cite{vine}, for each fixed $H \ge 1$, we can consider the discretized trajectory $\{(B_i(h), D_i(h))\}_{ h \in \Xi_H}$ of the $i$th feature, where $\Xi_H:=\{0, 1/H, \dots, (H - 1)/H\}$ describes a subdivision into slices of distance $1/H$. 

Now, we associate with each vine some statistic $\xi\big(\{(B_i(h), D_i(h))\}_{h \in \Xi_H}\big)$ and form a global test statistic by aggregating over all features in the window $Q_n$:
\begin{align}
	\label{long_eq}
	T_n = T(\XX_n) := \sum_i \xi\big(\{(B_i(h), D_i(h))\}_h\big).
\end{align}
For instance, the score $\xi$ could be the average number of slices in $\Xi_H$ where the feature is present, or the life time averaged over those features. In Section \ref{testdef_subsec}, we define precisely test statistics of this flavor.

%
%
\bee[Cross-sectional statistics]
\label{ex:cross}
Longitudinal statistics have the appeal of capturing the evolution of individual features over time. In practice, however, the trajectories might be observed only over a very limited number of slices so that vines  may extend only of a handful of slices. Hence,  computing general longitudinal statistics could be challenging. Nevertheless, by carefully choosing the score function $\xi$, we construct \emph{cross-sectional test statistics} that can still be evaluated seamlessly in such cases.

More precisely, for any measurable function $\xi'\co [0, \t]^2 \to [0, \ff)$, we consider 
$$\xi\big(\{(B_i(h), D_i(h))\}_{h \in \Xi_H} \big) :=\f1H \sum_{h \in \Xi_H}\xi'(B_i(h), D_i(h)),$$
implicitly setting $\xi'(B_i(h), D_i(h)) = 0$ if the $i$th feature is not present at height $h$. Then, forming the test statistic $T_n$ as in \eqref{long_eq}, and exchanging sums, we arrive at 
		\begin{align}
			\label{cross_eq}
			T_n = \f1H \sum_{h \in \Xi_H} \sum_i \xi'(B_i(h), D_i(h)) .
		\end{align}
		Here, we stress that to evaluate $T_n$ it is not necessary to track features, since the inner sum depends only  on the persistence diagram at time $h$. In fact, we can consider this estimator also as an integration with respect to the empirical measure associated with the persistence diagrams in the individual slices.

To make this precise, we may think of the persistence diagrams in a window as an empirical measures on $ [0, \t]^2$ averaged over the slices, i.e.,
	$$	\PD^q_n(\d b, \d d) := \f1H\sum_{h \in \Xi_H}\sum_{i\ge 1} \de_{(B_i^q(h),D_i^q(h))}(b, d).$$
Then, the test statistic \eqref{cross_eq} becomes 
	\begin{align}
		\label{cross_eq2}
		T_n :=	\lan \xi',  \PD^q_n\ran = \f1H \sum_{h \in \Xi_H} \sum_{i \ge 1} \xi'( B_i(h), D_i(h)).
	\end{align}
	For instance, taking 	$\xi'(b, d) := d - b$ yields total persistence averaged over the slices. As another example, taking the indicator 
$\xi'( b', d') := \one_{[0, b] \times [d, \t]}(b', d')$
	yields the \emph{slice-averaged $M$-bounded persistent Betti numbers}
$\b_n^{b, d} := \PD_n^q( [0, b], [d, \t])$.
\ene

\subsection{Stabilization conditions}
\label{cond_sec}
%
%
To develop goodness-of-fit tests based on the longitudinal test statistics described above, we show that they are asymptotically normal in large domains. On a very general level, proofs of asymptotic normality are typically based on certain dependence assumptions, and we now elaborate on what form these assumptions take in the present setting.

%
%
Although we fix $p = 3$, we keep the abstract notation, writing $\R^p$, $\R^{p'}$ instead of $\R^3$ and $\R^2$. This makes the arguments more transparent, and simultaneously opens the door towards future generalizations. 

%
%
We assume that $\XX_n$ emerges from a homogeneous Poisson point process $\PP= \{P_i\}$ with intensity $\lambda >0$ in the background through a construction rule that is stabilizing in the spirit of \cite{yukCLT}. This means that changing $\PP_n := \PP \cap Q_n$ far away from a position $x \in Q_n$ does not change the trajectories in the vicinity of $x$. More precisely, we assume that the family $\XX_n$ of trajectories is expressed as 
\begin{equation}\label{eq:XX}
\XX_n  = \{\TT(P_i, \PP_n)\}_{P_i \in \PP_{n - \sqrt n}}\end{equation}
	for some measurable and translation-covariant transformation $\TT$. In Section \ref{vor_sec} below, we present a detailed example for $\TT$ based on the Voronoi tessellation. Here, to control edge effects, we restrict to trajectories associated with Poisson points in the slightly eroded window $Q_{n - \sqrt n}$.

%
%
Moreover, we assume exponential stabilization in the vein of \cite{yogeshCLT, yukCLT}. This involves studying the effects of changing a point pattern $\vp$ close to a given location $x \in Q_n'$. More precisely, for any $r>0$ we put $Q_r(x):= (x + Q_r') \ti [0, 1]$, $Q(x) := Q_1(x)$ and for any finite $\AA \su Q_n \sm Q_r(x)$, $r<n$, we let 
$$\vp_{r, \AA} :=  \AA \cup\bi(\vp \cap Q_r(x) \bi)$$
be the configuration obtained from $\vp$ by placing $\AA$ outside of $Q_r(x)$, and for $\BB \su Q(x)$ we let 
$$\vp_{r, \AA, \BB} :=  \AA \cup \BB \cup\bi(\vp \cap (Q_r(x) \sm Q(x))\bi)$$
be the configuration where additionally $\BB$ is placed in $Q(x)$. Moreover, for two trajectories $f,\,g\in C([0,\,1],\, Q_n)$, we define their distance through
$
	\ms{dist}(f,\,g):=\inf_{h \le 1}\|f(h)-g(h)\|.
	$

	%
	%
\begin{definition}[Stabilization] A finite $\vp \su Q_n$ \emph{stabilizes for $x \in \R^{p'}$ at $r \ge 1$ relative to $Q_n$} if 
\begin{enumerate}[label=(\textbf{S\arabic*}),ref=(\textbf{S\arabic*})]
\item\label{S1}
We have 
$$\TT(P, \vp_{r, \AA, \BB}) = \TT(P, \vp_{r, \AA, \es})$$
 for any finite $\AA \su Q_n \sm Q_r(x)$, $\BB \su Q(x)$ and $P \in \vp_{r, \AA} \sm Q_r(x)$. Loosely speaking, for the configuration outside $ Q_r(x)$, what happens inside $ Q(x)$ is not relevant.
\item \label{S2}We have
$$\ms{dist}\big(\TT(P, \vp_{2r, \AA, \BB}), \TT(P', \vp_{2r, \AA, \BB})\big) > M$$
for any finite $\AA \su Q_n \sm Q_{2r}(x)$, $\BB \su Q(x)$, $P \in \AA$ and $P' \in \vp_{2r, \AA, \BB} \cap  Q_r(x)$. Loosely speaking, $M$-bounded features with points $P\in\AA$ cannot contain points inside $Q_r(x)$.
\item\label{S3} We have $\TT(P, \vp_{4r, \AA}) = \TT(P, \vp_{4r, \es})$
	for any finite $\AA \su Q_n\sm Q_{4r}(x)$, and $P \in \vp \cap Q_{2r}(x)$. Loosely speaking, trajectories in $Q_{2r}(x)$ will not be affected by the change to the configuration $\AA$. Here, we changed from $2r$ to $4r$ in order to harmonize well with condition \ref{S2}.
\end{enumerate}
	Formally, we also set $Q_\ff := \R^{p'}\ti [0, 1]$, and impose in this case additionally that the entire trajectories from~\ref{S3} be contained in $Q_{4r}(x)$. That is,
\begin{enumerate}[label=({\bf S\arabic*'}),ref=({\bf S\arabic*'})]
\setcounter{enumi}{2}
 \item \label{S3'} We have $\TT(P, \vp_{4r, \AA}) \su Q_{4r}(x)$
for any finite $\AA \su \R^p\sm Q_{4r}(x)$, and $P \in \vp \cap Q_{2r}(x)$.
\end{enumerate}

\end{definition}
Then, we let $R'(x, n; \vp)$ denote the smallest integer $r\ge 1$ for which $\vp$ stabilizes at $x$ relative to $Q_n$, and $R'(x, \ff; \vp)$ the smallest integer $r\ge 1$ for which $\vp$ stabilizes relative to $\R^p$. We define the \emph{stabilization radius} $R(x, n; \vp) := R'(x, n; \vp) \vee R'(x, \ff; \vp)$. {Note that while $R'(x, n; \vp)\le 2n$, it may happen that $R'(x, \ff;\vp) = \ff$}. In particular, we write $R(x, n):= R(x, n;  \PP)$ for the stabilization radius associated with the input point process. Finally, we need the notion of exponential stabilization.

%
%
\begin{definition}[Exponential stabilization]$\TT$ is \emph{exponentially stabilizing} if for some $c > 0$ one has
\begin{align*}
	\limsup_{r \to \ff} \sup_{\substack{n\ge 1 \\ x \in Q_n}}\f{\log(\P(R(x, n) > r))}{r^c} < 0.
\end{align*}
\end{definition}

\subsection{Asymptotic normality}
\label{norm_sec}

%
%
After having specified the conditions, we now state the asymptotic normality of scalar test statistics.

%
%
\bet[Asymptotic normality; scalar level]
\label{clt_thm}
Assume that $\TT$ is exponentially stabilizing and that the score function $\xi$ is bounded. Then, the test statistic 
$$\f{T_n - \E[T_n]}{n^{p'/2}}$$
converges in distribution to a normal random variable.
\ent

%
%
When dealing with cross-sectional statistics, previous work from \cite{svane}  already suggests that it should not only be possible to derive asymptotic normality separately for each scalar statistics $F$ but obtain a functional statement. Our second main result, Theorem \ref{fclt_thm}, will show that this is indeed the case. The key advantage of functional asymptotic normality of the persistent Betti numbers is that one can consider Kolmogorov-Smirnov-type statistics depending on a continuous range of birth/death-parameters.

To prove the functional CLT, we will need to impose further moment conditions. More precisely, we assume that the factorial moment densities $\rh(x_0, \dots, x_q)$ of the points $\{X_{i, n}(h)\}_{i \ge 1}$ in the slice $h \in \Xi_H$  are bounded in the sense that for every $q \ge 1$,
\begin{align}
	\label{mom_eq}
	\sup_{\substack{h \in \Xi_H\\ n \ge 1}} \sup_{x_0, \dots, x_q \in \R^p}\rh(x_0, \dots, x_q) \le C_{\r, q}\tag{{\bf M}} 
\end{align}
for some $C_{\r, q} > 0$.  Here, the factorial moment densities are determined by the disintegration property
\begin{equation}\label{eq:fac_mom}\E\Big[\sum_{\substack{X_{i_0, n},\dots, X_{i_q, n}  \in \XX_n}}f(X_{i_0, n}(h), \dots, X_{i_q, n}(h))  \Big] = \int_{\R^{dq}} f(x_0, \dots, x_q) \rh(x_0, \dots, x_q) \d(x_0, \dots, x_q)\end{equation}
	where $X_{i_0, n},\dots, X_{i_q, n}$ must be pairwise distinct, see \cite[Section 4.2]{poisBook}.

%
%
\bet[Asymptotic normality of $M$-bounded persistent Betti numbers; functional level]
\label{fclt_thm}
Assume that $\TT$ is exponentially stabilizing and that condition \eqref{mom_eq} holds. Then, as a function on $[0, \t]^2$, in the Skorohod topology, the recentered and rescaled $M$-bounded persistent Betti numbers 
$$\f{\bn - \E[\bn]}{n^{p'/2}}$$
converge in distribution to a centered Gaussian process.
\ent

%
%
\section{Example}
\label{vor_sec}
In this section, we present a specific model for the trajectories $\XX_n = \{X_{i, n}\}_{i \ge 1}$, namely as perturbed centroids of the 2D slices of a 3D Poisson Voronoi tessellation. More precisely, let $\PP_n = \{P_{i, n}\} \su Q_n$ be a homogeneous Poisson point process and consider the associated 3D Voronoi tessellation $\{\CC_{i, n}\}_{i}$ on $\PP_n$. Then, for each $h \in [0, 1]$ and each bounded cell $\CC_{i, n}$, we let $G_{i, n}(h)$ be the centroid of the section $\CC_{i, n} \cap (\{h\} \times \R^2)$ provided that the latter is non-empty. {We denote the entire collection of centroids at level $h$ as $G_n(h)$.}

We assume that, from the data, we only extract a slight perturbation of $G_{i, n}$. This added randomness  is justified from measurement errors and is also technically useful in the proofs. Hence, we set 
$$X_{i, n}(h) := G_{i, n}(h) + N_{i, n}(h),$$
where the $\{N_{i, n}(h)\}_{i, h}$ are assumed to be noise variables that are independent for different $i$ but may be correlated in $h$. For instance, we take $N_{i, n}(h)\sim \UUU(B(o, {\eta_0}))$ to be uniformly distributed in a ball of a small radius $\eta_0 < 1/4$. Here, we use the notation $B(x,\,m)\subset\R^d$ for the Euclidean ball around $x\in\R^p$ of radius $m>0$. We now verify condition \eqref{mom_eq} and exponential stabilization. Henceforth, to ease notation, $|A|$ denotes either the cardinality of a discrete set $A$ or the Lebesgue measure of $A$ if $A$ is a measurable subset of the Euclidean space.

{\em Condition \eqref{mom_eq}.}
Let $A_0, \dots, A_q \su Q_n$ be Borel sets whose diameter and volume we may assume to be at most 1. Then, we claim that 
\begin{align}
\label{vor_sec_eq}
	\E\Big[\Big|\big\{i_0,\dots, i_q \ge 1\text{ pairwise~distinct}\co X_{i_0, n}(h)\in A_0, \dots,X_{i_q, n}(h)\in A_q \big\}\Big|\Big]   \le C_q \prod_{j \le k}|A_j|,
\end{align}
with the constant $C_q$ depending neither  on $h$ nor $n$.

To prove this claim, we observe that since the support of the noise variables is of diameter at most $\et_0 \le 1/4$, we can have $X_{i_j, n}(h) \in A_j$ only if $G_{i_j, n}(h)$ is contained in the 1-neighborhood 
$$A_j^+ := \big\{ x + y\co  x \in A_j,\, y\in Q_{1/2}\big\}.$$
of $A_j$.  Since the noise variables are uniformly distributed, there exists some $c_0 > 0$ such that conditioned on the locations $\{P_i\}_{i \ge 1}$,
\begin{align*}
	\P\big(X_{i_0, n}(h)\in A_1, \dots,X_{i_q, n}(h)\in A_q\big\ba \{P_i\}_i\big) &\le c_0^k \one_{\{G_{i_1, n}(h)\in A_1^+, \dots, G_{i_q, n}(h)\in A_q^+\}}\prod_{j \le q}|A_j|.
\end{align*}
Hence, setting $A^+ := A_1^+\times\cdots\times A_q^+$ and  plugging in the right-hand side back into \eqref{vor_sec_eq}, it suffices to show that 
\mbox{$\sup_{ h \in [0, 1]}\E\big[ |G_n(h)^{q + 1} \cap A^+|\big] < \ff.$}
Moreover, by the H\"older inequality
, this further reduces to
$\sup_{\substack{ j \le k,\, h \le 1}}\E\big[|G_n(h) \cap A_j^+|^k\big] < \ff.$
	To prove this claim, fix some point $x \in A_j$ {and $h\in [0, 1]$}, and note that if $A_j$ is hit by $\ell \ge 1$ Voronoi cells, then there exists some $2 \le m \le 3n$ such that i) there are no points in $B(x, (m - 4)_+)$ and ii) there are at least $\ell$ Poisson points
in	$B(x,\,m) \sm B(x,\,{m - 2})$ (Lemma~\ref{lem:cells}).  Thus, for any $\ell \ge 2^8$,
\begin{align*}
	\P\big(|G_n(h) \cap A_j^+| = \ell\big) \le& \sum_{2 \le m \le 3n}\P\big(\PP_n \bi(B(x,\,m) \sm B(x,\,m-2)\bi) \ge \ell, \PP_n \cap B(x,\,(m-4)_+)= \es\big) \\
	\le&  \sum_{m \le \ell^{1/4}}\P\bi(\PP_n \bi(B(x,\,m) \sm B(x,\,m-2)\bi) \ge \ell\big) \\
	&+{\one\{3n \ge \ell^{1/4}\}}\sum_{\ell^{1/4} < m \le 3n}\P\bi( \PP_n \cap B(x,\,m-4)= \es\big),	
\end{align*}
where $\PP_n(A) := |\PP_n \cap A|$ denotes the number of points in a set $A$. 
Now, we bound the two summands separately. Since $\ell$ exceeds the expected number of Poisson points in $B(x,\,m) \sm B(x,\,m-2)$, the first summand decays exponentially fast in $\ell$ by the Poisson concentration inequality. The probability in the second sum is at most $\exp(-c_1 m^3 )$ for some universal constant $c_1>0$, so that the sum decays exponentially fast in $\ell^{3/4}$.

{\em Exponential stabilization.} We argue that exponential stabilization holds for the centroids, noting that a minor modification also gives the property after the perturbation.

To achieve this goal, let 
$$E_{r, x, n} := \Big\{Q_{\sqrt r}(\sqrt r z) \cap \PP \ne \es\co \text{for all }\sqrt r z \in \sqrt r \Z^{p - 1} \cap Q_n \cap Q_r(x)\sm Q_{r - 8 \sqrt r}(x)
\Big\}$$ 
denote the event that for each point $\sqrt r z \in \sqrt r \Z^{p - 1} \cap Q_n \cap Q_r(x)\sm Q_{r - 8 \sqrt r}(x)$ the $\sqrt r$-box $Q_{\sqrt r}(\sqrt r z)$ contains at least one point from $\PP$. We claim that 
\begin{equation}
\label{eq:union}
	\big\{E_{r/2, x, n} \cap E_{r, x, n} \cap E_{3r/2, x, n} \cap E_{5r/2, x, n}\big\} \su \big\{R'(x, n;\PP)\le r\big\},
\end{equation}
and similarly with $n$ replaced by $\ff$. Once we have proven this claim, the desired exponential stabilization will follow from a union bound in~\eqref{eq:union} and Lemma~\ref{lem:er} which gives, for two positive constants $c_1,\,c_2$
\begin{align*}
	&\P(E_{r,\,x,\,s}^c)\le c_1 r^{p'/2}\exp\big(-c_2 r^{p'/2}\big).
\end{align*}
To prove \eqref{eq:union}, we proceed in two steps: for all $s>0$
\begin{enumerate}[label=(\arabic*),ref=(\arabic*)]
 \item\label{claim1} under the event $E_{s, x, n}$, every cell centered in $Q_{s - 8 \sqrt s}(x)$ is contained in $Q_{s - 2 \sqrt s}(x)$;
 \item\label{claim2} under the event $E_{s, x, n}$, every cell centered outside $Q_{s + 2\sqrt s}(x)$ does not intersect $Q_{s - 2\sqrt s}(x)$. 
 \end{enumerate}
\definecolor{xdxdff}{rgb}{0.49019607843137253,0.49019607843137253,1.}
\definecolor{qqqqff}{rgb}{0.,0.,1.}
\definecolor{zzttqq}{rgb}{0.6,0.2,0.}
\begin{wrapfigure}{r}{0.4\textwidth}
\begin{tikzpicture}[scale=0.4,line cap=round,line join=round,>=triangle 45,x=1.0cm,y=1.0cm]
\clip(-3.4,-4.32) rectangle (9.68,5.5);
\fill[color=zzttqq,fill=zzttqq,fill opacity=0.1] (-1.58,5.44) -- (-1.58,-3.56) -- (7.42,-3.56) -- (7.42,5.44) -- cycle;
\fill[color=zzttqq,fill=zzttqq,fill opacity=0.1] (1.,3.) -- (1.,-1.) -- (5.,-1.) -- (5.,3.) -- cycle;
\draw [color=zzttqq] (-1.58,5.44)-- (-1.58,-3.56);
\draw [color=zzttqq] (-1.58,-3.56)-- (7.42,-3.56);
\draw [color=zzttqq] (7.42,-3.56)-- (7.42,5.44);
\draw [color=zzttqq] (7.42,5.44)-- (-1.58,5.44);
\draw (9.,3.)-- (3.,0.);
\draw [color=zzttqq] (1.,3.)-- (1.,-1.);
\draw [color=zzttqq] (1.,-1.)-- (5.,-1.);
\draw [color=zzttqq] (5.,-1.)-- (5.,3.);
\draw [color=zzttqq] (5.,3.)-- (1.,3.);
\draw [dash pattern=on 5pt off 5pt] (3.,3.)-- (9.,3.);
\draw [dash pattern=on 5pt off 5pt] (3.,3.)-- (3.,0.);
\draw [->] (-1.58,1.04) -- (1.,1.04);
\draw [->] (1.,1.04) -- (-1.58,1.04);
\draw (-1.08,2.4) node[anchor=north west] {$3\sqrt s$};
\draw [->] (-1.9,5.48) -- (-1.9,-3.48);
\draw [->] (-1.9,-3.48) -- (-1.9,5.48);
\draw (-3.2,0.5) node[anchor=north west,rotate=90] {$s-2\sqrt s$};
\draw [->] (1.02,-1.5) -- (5.06,-1.5);
\draw [->] (5.06,-1.5) -- (1.02,-1.5);
\draw (1.4,-1.7) node[anchor=north west] {$s-8\sqrt s$};
\begin{scriptsize}
\draw [fill=qqqqff] (9.,3.) circle (2.5pt);
\draw[color=qqqqff] (9.2,3.36) node {$P'$};
\draw [fill=xdxdff] (7.42,3.) circle (2.5pt);
\draw[color=xdxdff] (7.68,3.36) node {$P''$};
\draw [fill=qqqqff] (3.,0.) circle (2.5pt);
\draw[color=qqqqff] (2.74,0.36) node {$P$};
\draw [fill=xdxdff] (3.,1.) circle (2.5pt);
\draw[color=xdxdff] (3.54,1.36) node {$X$};
\end{scriptsize}
\end{tikzpicture}
\caption{Sketch for the proof of claim~\ref{claim1}}
\label{proof1_drawing}
\end{wrapfigure}
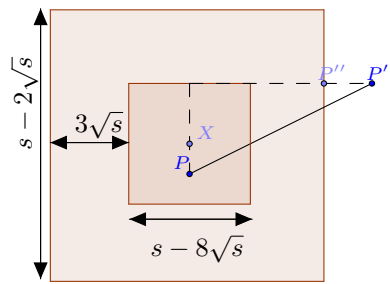
In particular, for sufficiently large $s > 0$,  under the event $E_{s/2, x, n}$ every cell centered in $Q_1(x)$ is contained in $Q_{s/2 - \sqrt{2s}}(x)$, and every cell centered outside $Q_{s/2 + \sqrt{2s}}(x)$ does not intersect $Q_{s/2 - \sqrt{2s}}(x)$ by~\ref{claim2}. Moreover, for sufficiently large $s > 0$, by~\ref{claim1} the event $E_{s, x, n}$ forces every cell centered inside $Q_{s/2 + \sqrt{2s}}(x)$ to be contained in $Q_{s - 2\sqrt s}(x)$. 
Therefore, under the event $E_{s/2, x, n} \cap E_{s, x, n}$, no cell centered in $Q_1(x)$ can intersect a cell with a point outside $Q_s(x)$, thereby establishing property~\ref{S1}. 
 Similar arguments show properties~\ref{S2},~\ref{S3} and \ref{S3'}.

We now provide the detailed derivations for claims~\ref{claim1} and~\ref{claim2}. For claim \ref{claim1},  let $C$ be a cell centered at $P \in Q_{s - 8 \sqrt s}(x)$, let $P'$ be an arbitrary point outside $Q_{s - 2\sqrt s}(x)$ and let $P''$ be its closest point on the boundary of $Q_{s - 2\sqrt s}(x)$. In particular, $|P' - P| \ge |P' - P''| + 3\sqrt s$. On the other hand, under the event $E_{s, x, n}$ the distance from $P'$ to one of the centers guaranteed by $E_{s, x, n}$ is at most $|P' - P''| + \sqrt {3 s} < |P' - P|$, so that $P' \not \in C$. The situation is sketched in Figure~\ref{proof1_drawing}.

For claim \ref{claim2},  let $C$ be a cell centered at $P \in \R^p \sm Q_{s + 2 \sqrt s}(x)$, let $P'$ be an arbitrary point inside $Q_{s - 2\sqrt s}(x)$ and let $P''$ be its closest point on the boundary of $Q_{s - 2\sqrt s}(x)$. Then, as before, $|P' - P| \ge |P' - P''| + 2\sqrt s$. On the other hand, under the event $E_{s, x, n}$, again the distance from $P'$ to one of the centers guaranteed by $E_{s, x, n}$ is at most $|P' - P''| + \sqrt {3 s} < |P' - P|$.

\section{Simulation study}
\label{simulation_sec}
In this section, we elucidate through a simulation study how to design a goodness-of-fit test based on the test statistics from Section \ref{stat_sec}. To that end, in Section \ref{set_sec}, we first describe the general set-up of this study including the considered null model, alternatives and test statistics. Then, Section \ref{expl_sec} provides an exploratory analysis illustrating the distribution of the proposed test statistics under the null model and the alternatives. Finally, Section \ref{power_sec} contains a more detailed and structured investigation concerning the test power.

%
%
\subsection{Simulation set-up}
\label{set_sec}
First, we describe the general set-up of the simulation study. That is, we present the null model, the alternatives as well as specific choices for test statistics. 

\subsubsection{Null model}
Recalling the setting of Example \ref{vor_sec}, we study a null model for the trajectories traced by the centroids of the 2D-slices of a 3D-Poisson Voronoi tessellation. For this simulation study, the null model is a Poisson Voronoi tessellation with cell centers given by a Poisson process with intensity $\lambda=2.18 \cdot 10^{-4}$ in a $170\times170\times85$-sampling window. The window height and the intensity of the point process are chosen to mimic the data set in Section \ref{data_sec}. Due to the asymptotic normality, it suffices to compute mean and variances under the null model in the chosen sampling window, whose extent in the $x$-$y$-directions is far smaller than that found in the real data set. After constructing the tessellation, we take 9 slices parallel to the $x$-$y$-plane with a fixed spacing of 4 between them.
The choice of the number of slices and of the spacing among them is guided by the exigency of mimicking a real serial slicing procedure, explained in more details in Section \ref{data_sec}.

\subsubsection{Alternatives}
In order to design a set of alternatives to the null model, we retain the basic assumption of taking multiple 2D slices of a 3D Voronoi tessellation but vary the set of cell generators. More precisely, we deviate from the Poisson model by considering two classes of point processes with interactions, namely a more regular and a more clustered point pattern. We refer the reader to \cite{daley2007II} for an overview of repulsive and clustered point processes.

As alternative to the null hypothesis of the Poisson process, we consider:
\beit
\im A \emph{Mat\'ern hard-core process}, which is based on a Poisson point process where any two points at distance smaller than a threshold $R > 0$ are removed.
\im A \emph{Mat\'ern cluster process}, where we first distribute a certain number $n_{\ms{cl}}$ of cluster centers in the window. Each cluster center generates offspring points according to a Poisson point process with intensity $\la_{\ms{cl}}$ inside a ball of radius $R$.
\enit

We work with two different parameter sets for the Mat\'ern hard-core and three different parameter sets for the Mat\'ern cluster point process, see Table \ref{tab:param}. The parameters of the different point processes are chosen so as to lead to the same expected number of points in the sampling window.

\begin{table}[h!]
\centering
\caption{Parameter choices for the alternatives.}
\label{tab:param}
\begin{tabular}{lllllll}
  \hline
	& $\ms{HC}_1 $&$ \ms{HC}_2 $&$ \ms{CL}_1$&$ \ms{CL}_2$&$ \ms{CL}_3 $ \\ 
  \hline
	$R$ &5.25 & 5.95 & 42.5  &  42.5 &  42.5 \\
	$n_{\ms{cl}}$ &/ & / & 10 & 5 & 4\\
	$\la_{\ms{cl}}$ &/ & / & 10 & 20 & 25
\end{tabular}
\end{table}
\vspace{-.4cm}

\begin{figure}[!h]
\centering
	\includegraphics[width=.7\textwidth]{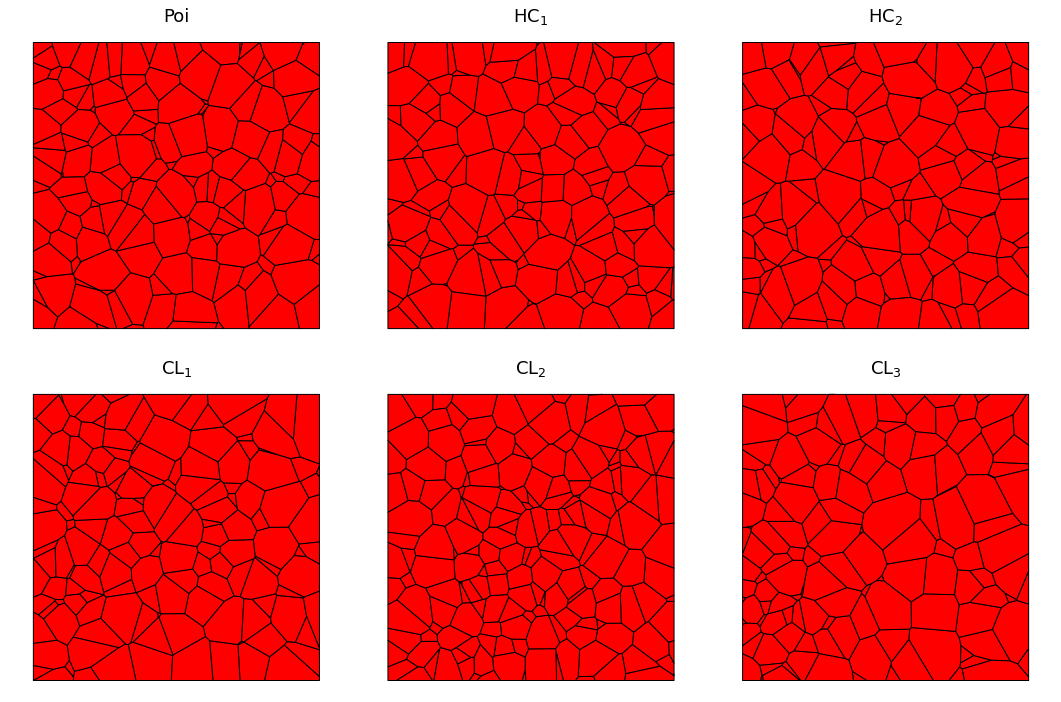}
\caption{2D slices of Voronoi tessellations generated from different point processes for the generator points.}
\label{fig:crossectioncomparison}
\end{figure}

Figure  \ref{fig:crossectioncomparison} illustrates the effects of varying the cell generators on a single 2D slice of the 3D Voronoi tessellation. We see that this variation induces rather subtle changes in the sliced cells, which look fairly similar under the different alternatives.

\begin{figure}[!h]
\centering
\includegraphics[width=\textwidth]{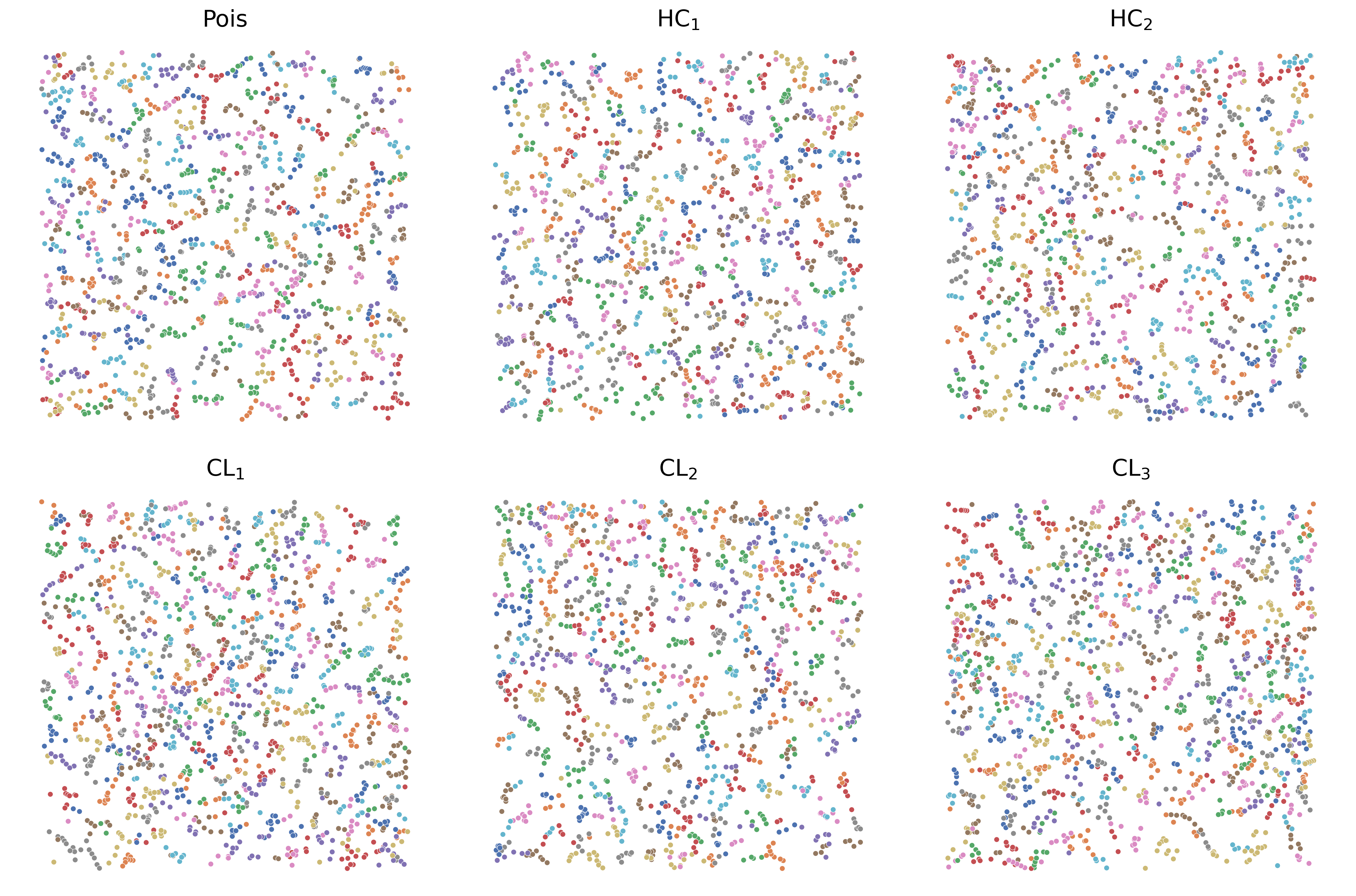}
\caption{$x$-$y$ coordinates of the centroids of the 2D sliced grains of the 3D Voronoi tessellations on the considered point processes. Different colors represent the different slices.}
\label{fig:centermasscomp}
\end{figure}

Furthermore, we compute the persistence diagrams of the different 2D slices using the centroids of the different slices  (shown in Figure \ref{fig:centermasscomp}) as input point cloud. This illustrates that although typically the center locations are rather close between the different slices, occasionally larger movements can be noted.

\subsubsection{Test statistics}
\label{testdef_subsec}
To test deviations of data from the null model, we rely on three different test statistics.

\been
\im {\bf Cross-sectional total persistence.} To begin with, we present an example of a cross-sectional statistics in the sense of Example \ref{ex:cross}. In words, we extract global information from each of the persistence diagrams computed for every slice and then aggregate these quantities into a single characteristic for the dataset. More specifically, we put
$$
\TTP^q:=\f1H\sum_{h \le H}\f1{|W|}\sum_{i\le n_h} \big(D_i^q (h) - B_i^q (h)\big).
$$
where $H$ is the total number of slices and $n_h$ the total number of features in slice $h$. That is, in the definition \eqref{cross_eq}, we take $\xi' (b, d) := d - b$ and then normalize the resulting statistic by the volume of the sampling window. %
\im {\bf Vine-based persistence.} Next, we discuss an example of a general longitudinal vine-statistic in the form of \eqref{long_eq}. Here, we average the life times of the features represented by the vine in the different slices over all slices where this vine is present. More precisely, we put 
$$\TM^q :=\f1{|W|}\sum_{i \le N^q}\f1{n_{h (i)}^q}\sum_{h\le H}  (D_i^q (h) - B_i^q (h)),$$
where  $N^q$ is the total number of unique features of dimension $q$ observed in all the slices $H$ and $n_{h (i)}^q$ is the number of slices in which the $i$th $q$-feature is visible. In order to represent $\TM^q$ in the general form of \eqref{long_eq}, we may choose $\xi\big (\{ (B_i^q (h), D_i^q (h))\}_h\big) := \f1{n_{h (i)}^q}\sum_{h\le H}  (D_i^q (h) - B_i^q (h))$. 
\im{\bf Ripley $K$-function.} Finally, in order to compare TDA-based test statistics to classical test statistics from spatial statistics, we discuss an example derived from Ripley's $K$-function.
More precisely, we let
$$
T_{\Rip}:=\int_0^{r_{\Rip}}\hat K_{\ms{pool}} (r)\d r,
$$
where, $\hat K_{\ms{pool}} (r)$ is the {pool}ed Ripley's $K$-function that combines the estimates of the Ripley's $K$-function in each considered slice. 
\enen

\subsection{Exploratory analysis}
\label{expl_sec}
\begin{figure}[!h]
	\centering
\includegraphics[width=\textwidth]{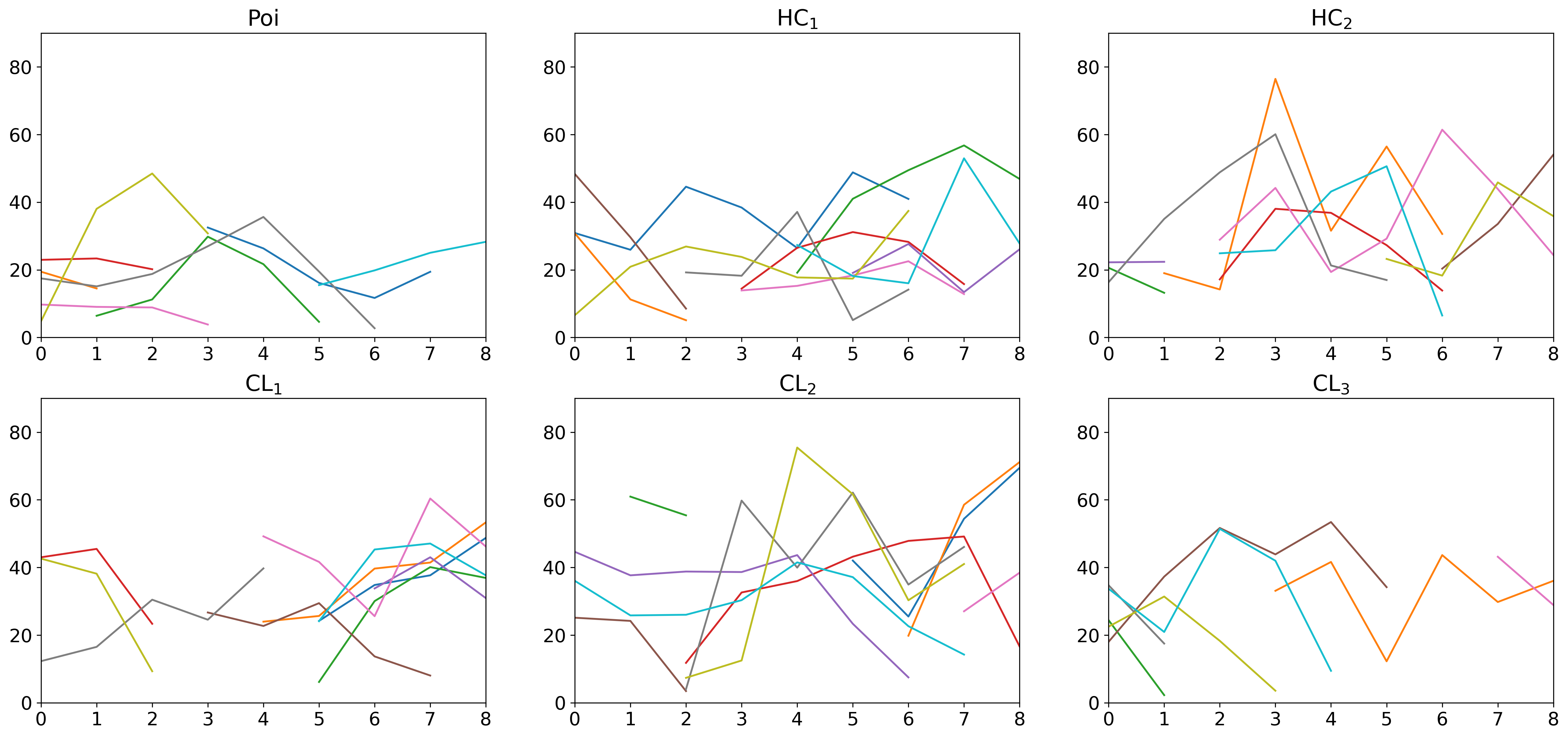}
	\caption{Samples of vines in dimension 0. Lines indicating the same persistence point observed in more than one slice.} 
\label{fig:pervincomp}
\end{figure}

\begin{figure}[!h]
	\centering
\includegraphics[width=\textwidth]{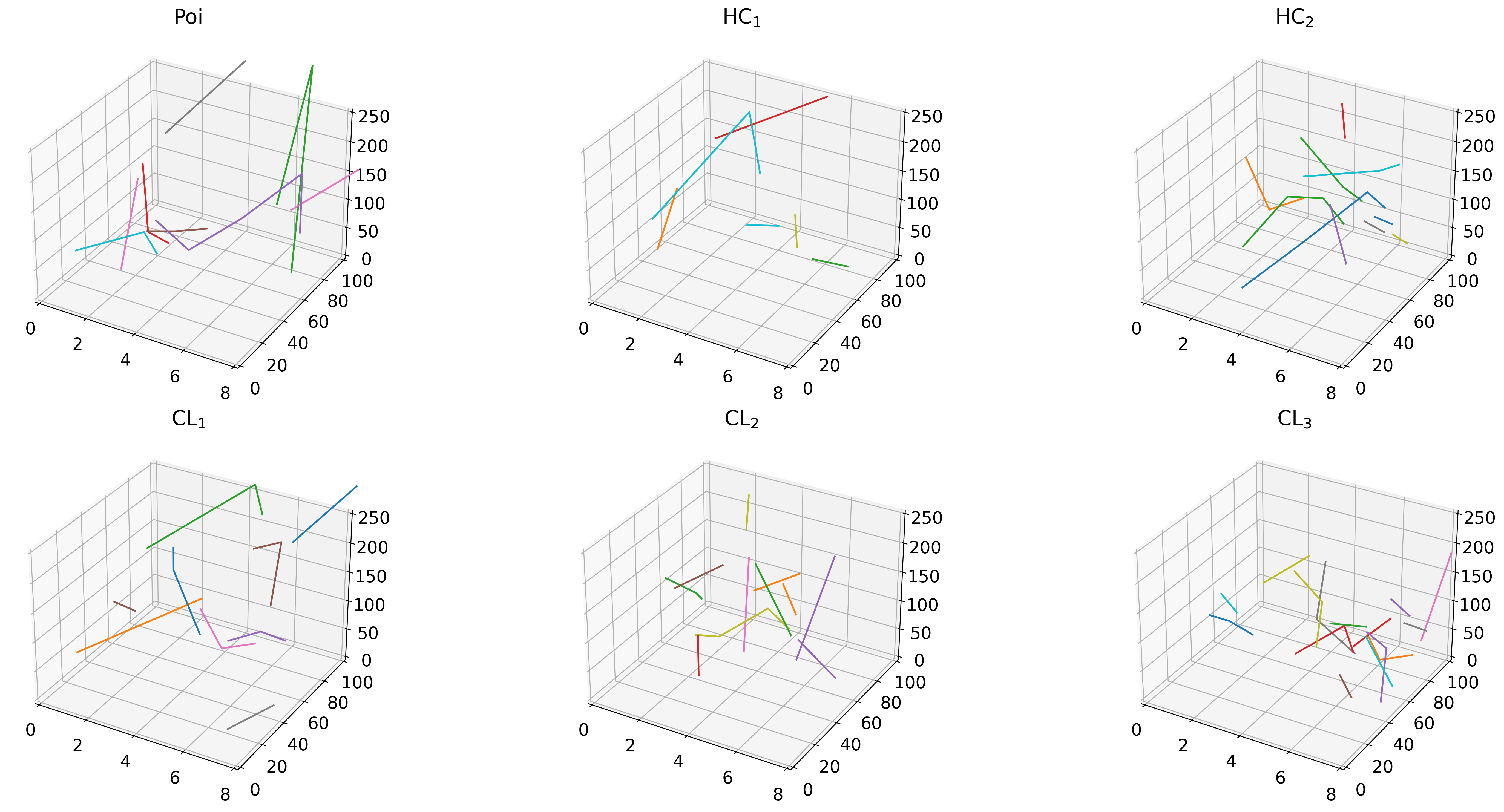}
	\caption{Samples of vines in dimension 1. Lines indicating the same persistence point observed in more than one slice.} 
\label{fig:pervincomp1}
\end{figure}

Now, we consider the setup of Example \ref{vor_sec}. First, Figures \ref{fig:pervincomp} and \ref{fig:pervincomp1}  respectively illustrate the 0-vines and 1-vines extracted from the persistence vineyard, obtained combining the persistence diagrams of each slice computed using the centroids shown in Figure \ref{fig:centermasscomp}. We see that both in the null model and in the alternatives the 0-vines occur in a variety of different lengths and present different trends. On the other hand, 1-vines are much shorter.

 This impression is also reinforced in Table \ref{tab:vinelength}, which shows the mean and the standard deviation of the vine length  ($L_v$) for the different Voronoi models. The lengths are rather similar between the alternatives but the average length of 0-vines exceeds that of 1-vines.
\begin{table}[h!]
\centering
	\caption{Left: Mean  (std.~dev.) of the {vine length} for the different Voronoi tessellations  (values based on 5,000 diagrams); Right: Mean  (std.~dev.) of the different test statistics in the Poisson-Voronoi model.}
\label{tab:vinelength}
\begin{tabular}{lcc}
	& {$\ms{dim}=0$}& {$\ms{dim}=1$}\\
  \hline
$\ms{PV}$ &2.771 (1.836)& 0.636 (0.878)\\
$\ms{HC}_1$ & 2.768 (1.824)& 0.643 (0.877)\\
$\ms{HC}_2$ &2.764 (1.821)& 0.643 (0.876)\\
$\ms{CL}_1$ &2.785 (1.851)& 0.635 (0.886)\\
$\ms{CL}_2$ & 2.775 (1.854)& 0.633 (0.888)\\
$\ms{CL}_3$ & 2.771 (1.856)& 0.630 (0.891)	
\end{tabular}\quad
	\begin{tabular}{lcc}
&  $170 \times 170 \times 85$ &  $140 \times 140 \times 140$ \\ 
  \hline
	$\TTP^0$ &1.618 (0.016) & 1.629(0.021)\\
	$\TTP^1$  & 0.899 (0.016) & 0.844(0.020)\\
	$\TM^0$ &0.397 (0.008) & 0.401(0.010)\\
	$\TM^1$ & 0.503 (0.017) & 0.474(0.020)\\
	$T_{\Rip}$  &355.507 (6.137) & 331.654(7.277) \\
	\phantom{aa}
\end{tabular}
\end{table}

\subsubsection{Mean and standard deviation under the null model}
The mean and variance of $\TTP^q$ and $\TM^q$ under the null model are computed using a simulation based on 5,000 Poisson-Voronoi tessellation for which 9 different slices with a spacing of 4 are considered. 

By Theorem \ref{clt_thm}, the statistics $\TTP^q$ and $\TM^q$ are asymptotically normal so that knowing the mean and variance allows us to construct a deviation test whose nominal confidence level is asymptotically exact. Results on the asymptotic normality related to Ripley's $K$-function in the single section case are obtained in \cite{ripley}.
We stress that our results on the asymptotic normality can be applied to a window of any fixed height as long as the size in the $x$- and $y$-directions are sufficiently large. On the other hand, it is also attractive to understand to what extent test statistics for windows of a given height generalize to other windows. This is not automatic since for fixed heights, the statistics may still be subject to finite-size effects. In Table \ref{tab:vinelength} (right), we compare the test statistics in the original $170 \times 170 \times 85$-window to the one computed in a cube of side length 140. We see that when moving to the cubical window, the mean-values for the statistics corresponding to features deviate by around 7\% or less from the values obtained in the original window. 

As a conclusion, we note that the approximation resulting from ignoring effects coming from a finite height is already accurate for moderately large sampling windows. We also note that in a different context, edge effects and the dependence on the number of observed cells in 2D was observed in the analysis in \cite{hahn,vittorietti}. There, this dependency is stated and the quantiles of the test statistics for different sample sizes  (e.g.~50, 100, 150 visible cells in 2D) are reported.

\begin{figure}[!h]
\centering
\includegraphics[width=\textwidth]{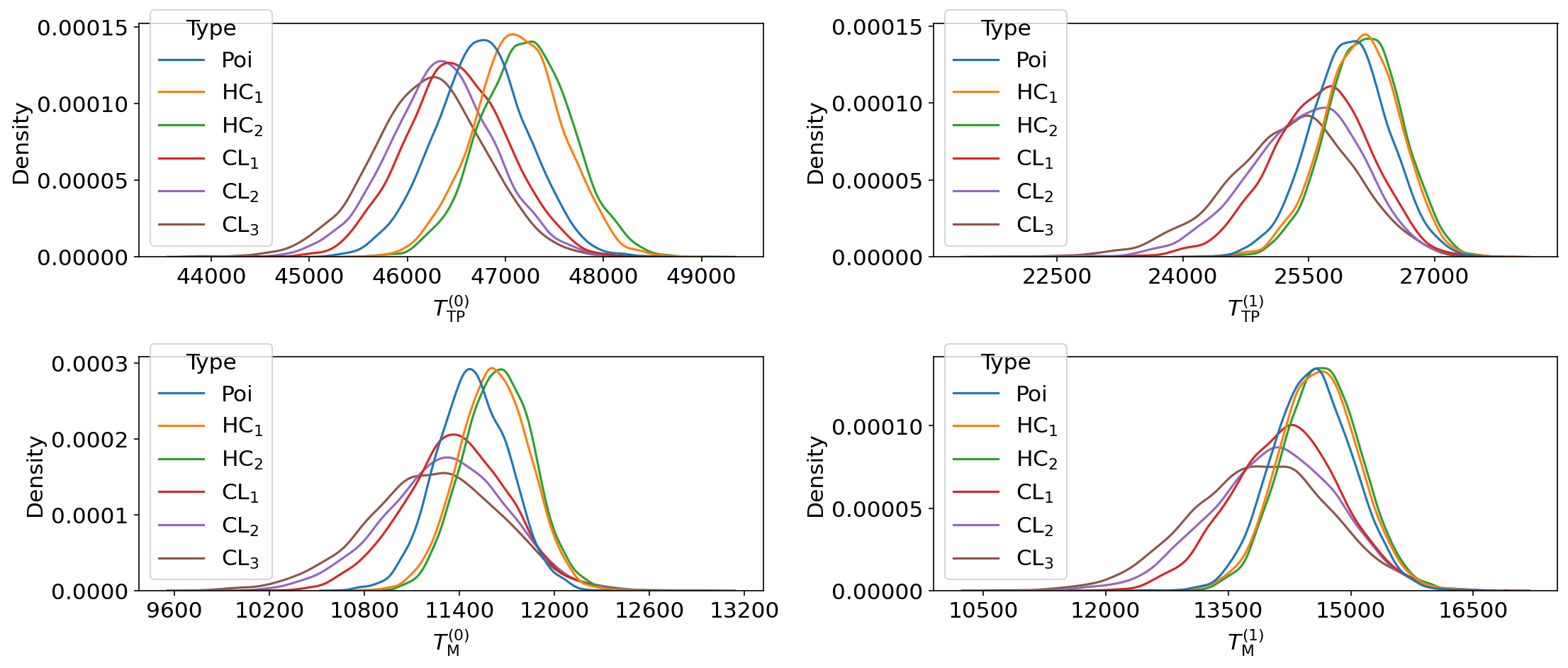}
	\caption{Distribution of the test statistics for the null model  (blue) and the alternatives based on 5,000 realizations under the null model and the alternatives. }
\label{fig:cdfcomp}
\end{figure}

Figure \ref{fig:cdfcomp} compares the empirical cumulative distribution function of the test statistics computed under the null and the alternative models.
For all of the considered test statistics, the distributions of the Mat\'ern cluster patterns differ clearly from the one of the Poisson null model. On the other hand, the situation for the Mat\'ern hard-core patterns is more subtle. While we can observe differences for features in degree 0, there is very little difference when considering features in degree 1.
%

\subsection{Power analysis}
\label{power_sec}
Before presenting the results, we explain
in detail how to run the test. First, we choose
a significance level of $\alpha=0.05$. Second, we generate 5,000 realizations from the null model and from the alternatives. The Voronoi tessellations under the alternatives are the ones described at the beginning of the section. 
To analyze the type I and II errors, we draw Table \ref{tab:rejectionrates} that contains the rejection
rates of this test setup. 

Under the null model the rejection rates are close to the nominal 5\%-level thereby illustrating that already for moderately large point
patterns, the approximation by the Gaussian limit is accurate. For the Mat\'ern hard-core point pattern, we see that the rejection rates for $\TTP^0$ are higher than those for $\TTP^1$, whereas for the Mat\'ern cluster point patterns, the situation is reversed. Moreover, the rejection rates for the $\TTP$-statistics are slightly higher than those for the $\TM$-statistics. We also observe that the pooled Ripley-statistics $T_{\Rip}$ is very powerful for detecting the clustered point patterns. On the other hand for the hard-core point patterns, the test statistic $\TTP^0$ has a bit higher rejection rate.
\begin{table}[!hbp]
\centering
\caption{Rejection rates for the multi-slice test statistics under the null model and the alternatives}
\label{tab:rejectionrates}
{\begin{tabular}{lrrrrrr}
T&   $\ms{PV}$ & $\ms{HC}_1$& $\ms{HC}_2$& $\ms{CL}_1$& $\ms{CL}_2$&$\ms{CL}3$  \\ 
  \hline
$\TTP^0$&	5.08\%& 11.96\%& 15.94\%& 10.70\%&16.58\%&24.50 \%\\
$\TTP^1$&5.30\%& 5.10\%& 5.64\%& 18.24\%&28.77\%&38.01 \%\\
	$\TM^0$&5.00\%& 8.30\%& 9.46\%& 18.84\%&27.47\%&34.85 \%\\
	$\TM^1$&4.74\%& 5.10\%& 5.02\%& 19.92\%&27.83\%&37.15 \%\\
	$T_{\Rip}$&4.94\%& 7.82\%& 9.76\%& 52.01\%&70.17\%&79.50\%\\
\end{tabular}}
\end{table}

\subsection{Extensions and variations}
In this section, we discuss three possible extensions and variations of the testing methodology described above. 

\subsubsection{Single-slice testing.}
To highlight the difference between testing based on multiple slices and testing based just on one slice, in Table \ref{tab:rejectionratessingolarsec} the rejection rates obtained using just one slice are reported: the test based on multiple slices presents higher rejection rates and it does not turn to be too conservative.
\begin{table}[!h]
\centering
\caption{Rejection rates for the single-slice test statistics under the null model and the alternatives.}
\label{tab:rejectionratessingolarsec}
{\begin{tabular}{lrrrrrr}
T&   $\ms{PV}$ & $\ms{HC}_1$& $\ms{HC}_2$& $\ms{CL}_1$& $\ms{CL}_2$&$\ms{CL}_3$  \\ 
  \hline
$\TTP^0$&4.70\%& 5.54\%& 6.88\%& 6.58\%&7.94\%&9.92 \%\\
$\TTP^1$&4.96\%& 5.02\%& 5.62\%& 9.20\%&12.90\%&17.16 \%\\
$T_{\Rip}$&5.38\%& 5.10\%& 5.10\%& 25.55\%&40.79\%&52.77\%\\
\end{tabular}}
\end{table}

\subsubsection{Cell vertices instead of centroids}
Moreover, in the present simulation study, we constructed the persistence diagrams based on the centroid of the cell slices. A concern of this methodology is that subtle differences between point patterns could be lost when moving from the sliced Voronoi cells to the centroids. Therefore, we also experimented with  computing the test statistics $\TTP^i$ when computing the persistence diagram on the basis of the cell vertices instead of the cell centroids. However, we found that this modification did not improve the testing power. We hypothesize that this is due to subtle dependencies induced by the slicing procedure. Even if for 3D cells the centroid may not be highly informative, this situation changes when working with the sequence of centroids obtained from the multiple 2D slices.

\subsubsection{Labeling algorithm}
In real data we have to deal with the problem of identifying to which grain the centroid belongs. Often, additional information on the orientation of the grains helps reconstructing the inner grain structures of the slices.
Alternatively, we propose a simple algorithm for the labels assignment. We start by computing the pairwise distance matrix between the set of centroids that belong to adjacent slices. We then assign the label of the points of the first slice to the points of the second slice if and only if the distance between two points does not exceed a specific threshold and it is minimum with respect to all the other points. 

To assess the effectiveness of this algorithm, in Table \ref{tab:lab}, we present the reconstruction error of the vines for the different models under consideration.  Although the simple algorithm does not lead to a perfect reconstruction, it still succeeds in recovering a substantial proportion of the vines across the models. Moreover,  in Table \ref{tab:lab} we also note a decrease in the testing power of $\TM^0$ and $\TM^1$ when using the reconstructed labels instead of the ground truth. This illustrates that the test statistics $\TM^0$ and $\TM^1$ are more useful in settings featuring a relatively densely arranged set of slices. Their power may deteriorate if the slices are so far apart that a precise vine reconstruction may no longer be feasible. Note that $\TTP^0$, $\TTP^1$ and $T_{\Rip}$ instead are not affected by the reconstruction algorithm.

\begin{table}[!h]
\centering
\caption{Reconstruction error $\ms{Err}_{\ms{Rec}}$ for the actual vines under the reconstruction algorithm.  Rejection rates for the test statistics based on the reconstructions the null model and the alternatives.}
\label{tab:lab}
{\begin{tabular}{lcccccc}
&   $\ms{PV}$ & $\ms{HC}_1$& $\ms{HC}_2$& $\ms{CL}_1$& $\ms{CL}_2$&$\ms{CL}3$  \\ 
  \hline
$	\ms{Err}_{\ms{Rec}}$&21.4\%& 20.6\%& 20.4\%& 21.6\%&21.8\%&21.9\%\\
$\TM^0$&5.14\%& 5.04\%& 5.04\%& 8.08\%&9.52\%&12.96\% \\
$\TM^1$&4.64\%& 4.12\%& 4.66\%& 14.36\%&19.96\%&27.07\%
\end{tabular}}
\end{table}

\section{Analysis of material data}
\label{data_sec}
In this section, the use of Poisson-Voronoi diagrams for representing the microstructure of extra low carbon strip steel is tested. By relying on the sequential slicing technique via automated mechanical polishing, 3D electron backscatter diffraction (EBSD) measurements of this material are carried out.

In this technique, consecutive steps of sample preparation and EBSD scanning are employed to obtain 2D EBSD slices. More details about the material and the technique can be found in \cite{galan}.

\subsection{Exploratory analysis}

\begin{figure}[!h]
\centering
\includegraphics[scale=0.35]{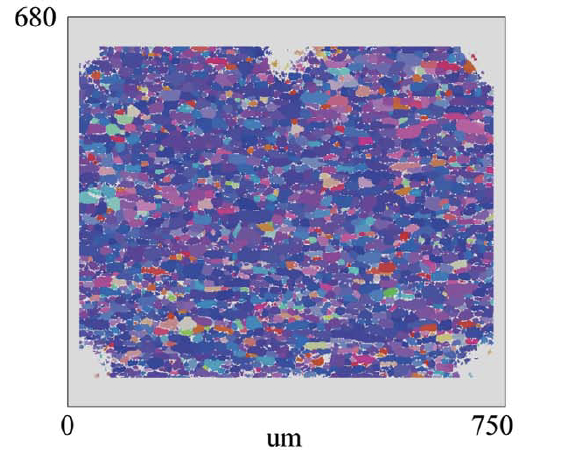}
\includegraphics[scale=0.5]{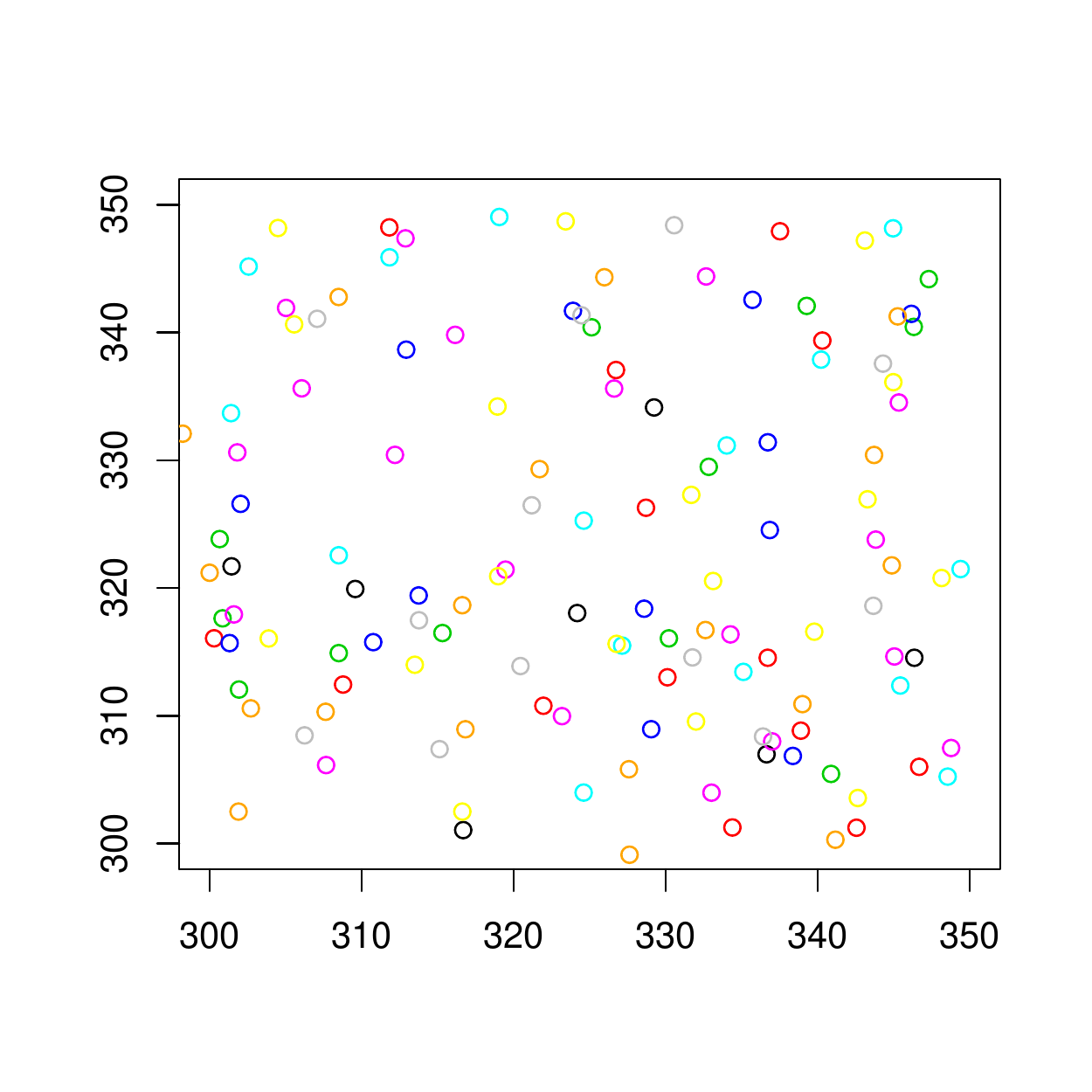}
	\caption{Experimental EBSD data (left), and $x$--$y$ coordinates of centroids of the 2D sliced cells of nine 2D EBSD scan (right). Different colors represent the different slices.}
\label{fig:ebsddata}
\end{figure}

\begin{figure}[!h]
\centering
\subfloat[]{\includegraphics[scale=0.35]{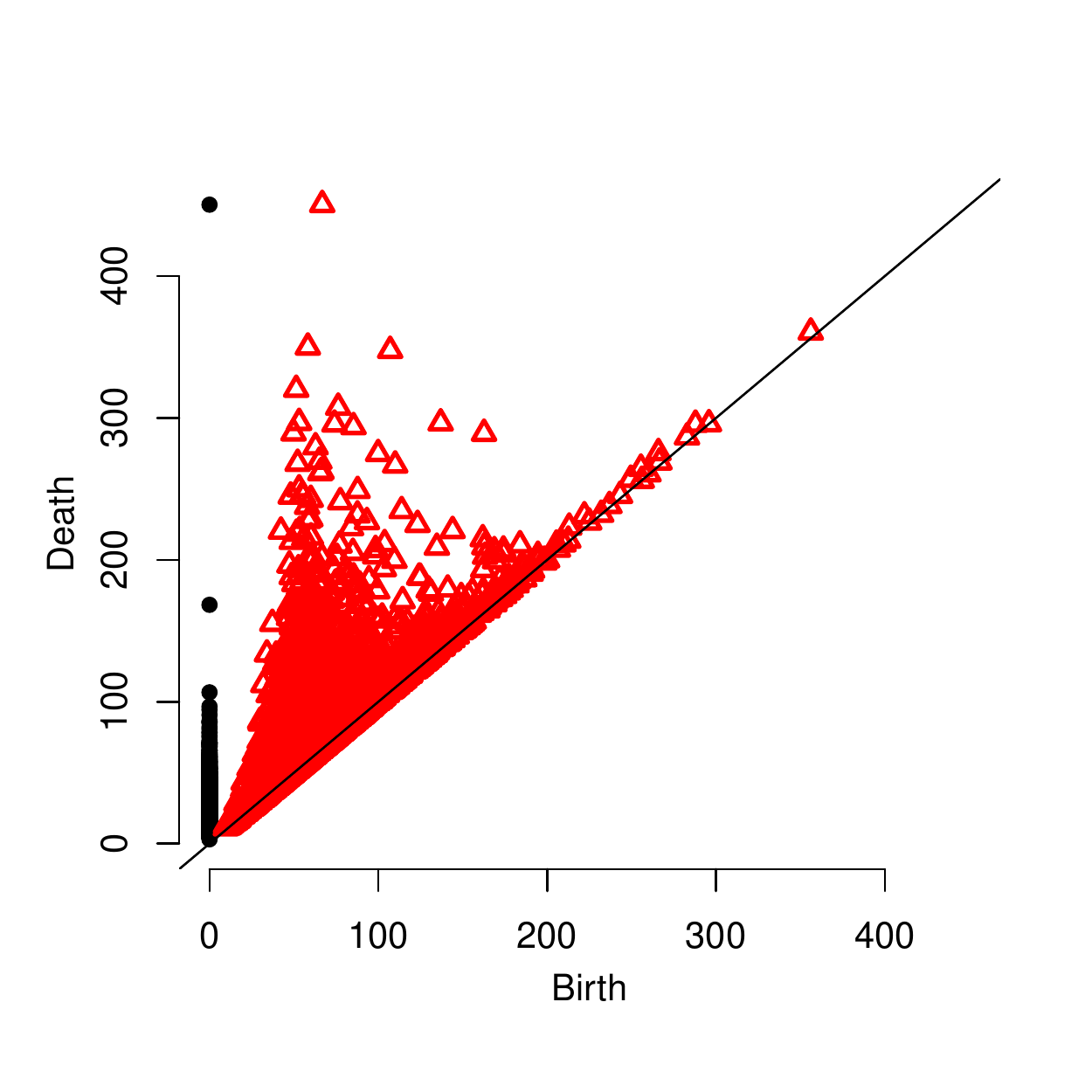}}
\subfloat[]{\includegraphics[scale=0.35]{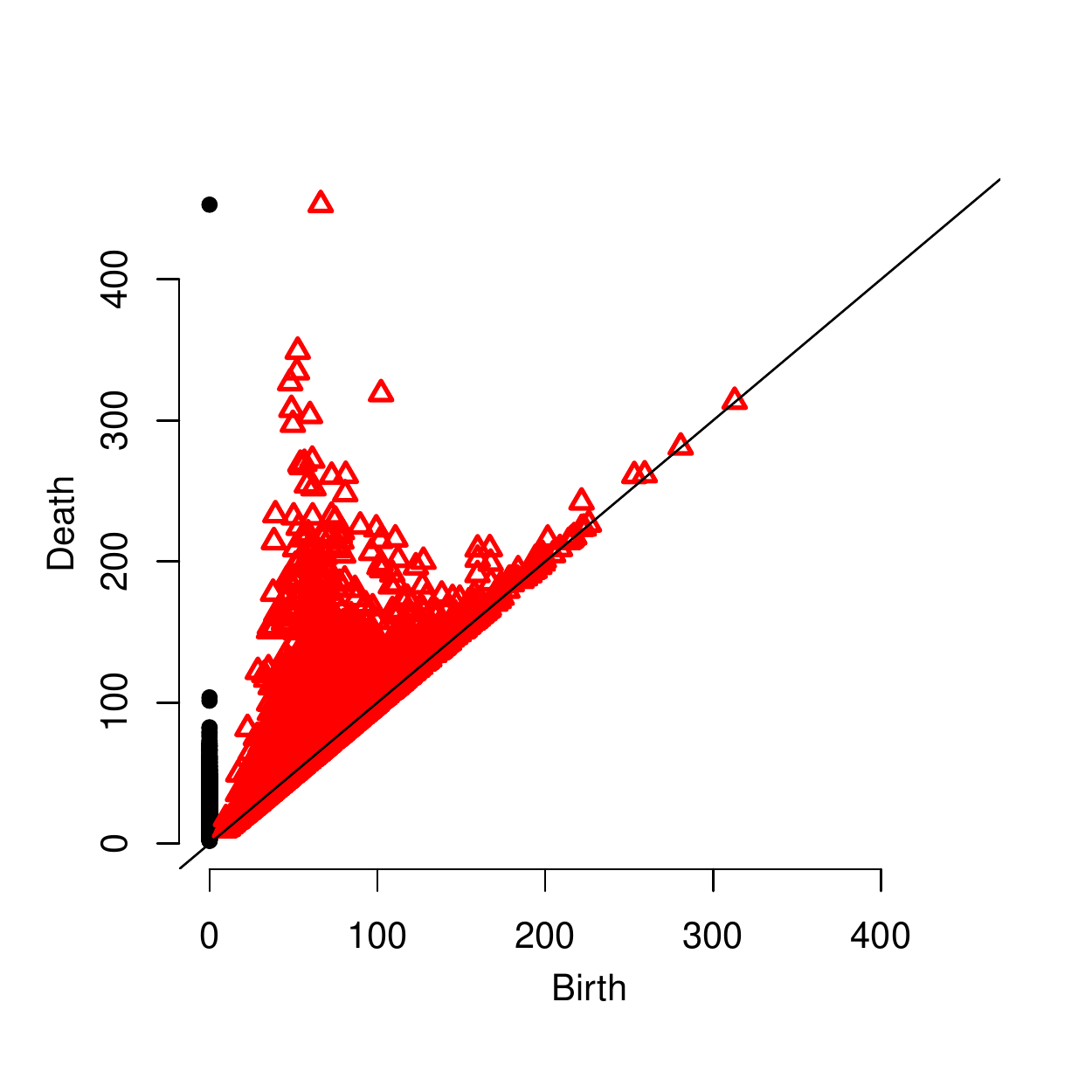}}
\subfloat[]{\includegraphics[scale=0.35]{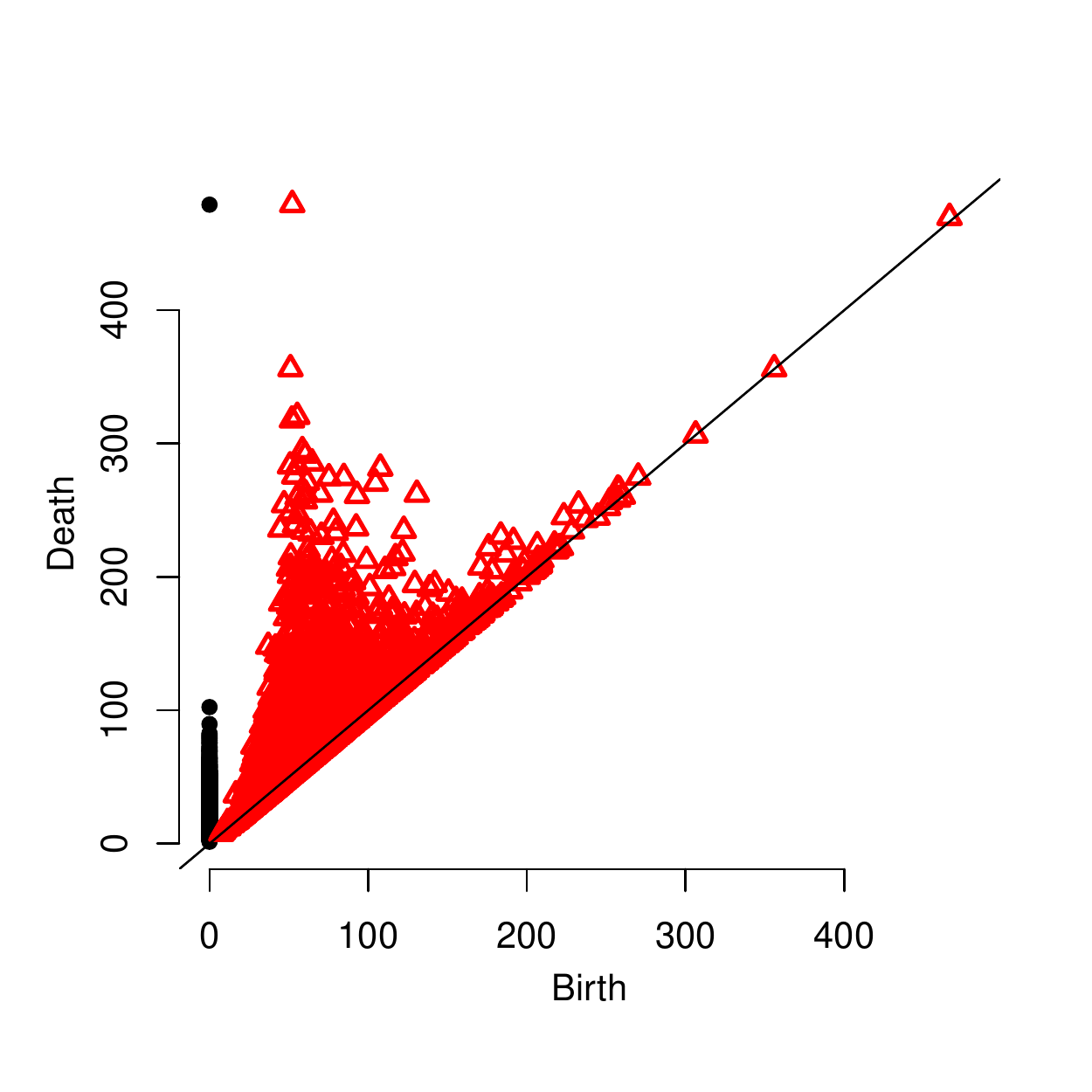}}
\caption{Persistence diagrams corresponding to three different 2D EBSD slices based on the point clouds of centroids of the 2D sectional cells.}
\label{fig:pdrealdata}
\end{figure}

The left panel in Figure \ref{fig:ebsddata} (taken by J.~G.~L\'opez \cite{galan}) illustrates a sample of the steel data consisting of 20 EBSD scans with a spacing of 4$\mu$m. On average $\approx 2,700$ grains ($sd\approx120$) are observed in each 2D EBSD scan.
In order to reduce edge effects, we proceed as in Section \ref{simulation_sec} and restrict the computation to nine 2D EBSD scans taken in the middle of the surface block.
In Figure \ref{fig:ebsddata} (right) the cell centroids  of the grains visible in nine 2D EBSD scans are plotted (a reduced area is considered for visualization purposes). Each set of centroids is then used for building a persistence diagram (Figure \ref{fig:pdrealdata}).

In total $24,008$ distinct features are observed in dimension 0 and $31,474$ in dimension 1.
The mean ``vine'' length is $0.021243$ and  $0.0004$,  highlighting that only few features are actually observed in more than one slice.
\newpage
\subsection{Poisson-Voronoi testing}
Under the impression of the previous visualizations, we now test the Poisson Voronoi hypothesis for the EBSD data. As in the previous section, we use the asymptotic normality of the test statistics under the null model. The test statistics are reported in Table \ref{tab:testrealdata}. 
\begin{table}[ht]
\centering
\begin{tabular}{lrrrrr}
	&  $\TTP^0$ &  $\TTP^1$ & $\TM^0$  & $\TM^1$&$T_{\ms{Rip}}$  \\ 
  \hline
	$z$-score&	23.00&   5.40 &94.75  & 30.05 & 81.8
\end{tabular}
\caption{$z$-scores associated with the different test statistics when the 2D EBSD slices are tested against a Poisson-Voronoi null model.}
\label{tab:testrealdata}
\end{table}

First, we note that all considered test statistics suggest a highly significant deviation of the data from the 3D null model of a Poisson-Voronoi tessellation. Second, we see that the test statistics $\TM^0$ and $\TM^1$ lead to $z$-scores that exceed by far the test scores corresponding to $\TTP^0$ and $\TTP^1$. This effect points to the potential that lies in the longitudinal statistics. Indeed, as observed at the beginning of the present section, on average, the vines in the data set are far shorter than the ones observed in the null model. Thus, the deviations from the null model are far more pronounced in the longitudinal test statistics $\TM^0$ and $\TM^1$ than they are in the cross-sectional test statistics $\TTP^0$ and $\TTP^1$. 

\section{Conclusion and perspective}
\label{conc_sec}
In this work, we applied the concept of persistence vineyards to develop goodness-of-fit tests for data given in the form of 2D slices of a 3D data set. These tests rely on the asymptotic normality of the considered test statistics, which we establish under suitable moment and stabilization conditions. The potential of the new methodology is illustrated through a simulation study and a dataset describing the microstructure of extra low carbon strip steel. We stress that the present work is to be considered only as the first step towards using the tools from TDA to develop statistical tests for slice-based datasets. We discuss now possible avenues for further research.  

First, motivated by the shape of the specific dataset, we established the asymptotic normality for domains that are largein the $x$- and $y$-direction but are fixed in the $z$-direction. For other datasets, it may be that all three dimensions are large, and it would be of interest to extend the proof of asymptotic normality for such settings. This would have the additional benefit that under a Poisson-Voronoi null model, it suffices to compute the mean and variance of the test statistics in Section \ref{simulation_sec} for a fixed point-process intensity. The values for other intensities can then be obtained through a suitable scaling by a power of the intensity. It is also interesting to extend the results to higher dimensions, power-law correlations and other stochastic model for representing microstructures.

Second, also our decision to choosing the centroids of the sliced cells as a basis for the persistence computation is motivated from the constraints found in the considered dataset. As illustrated in Figure \ref{fig:ebsddata}, the  cells in the material data can be shaped rather irregularly. Since the cell centroids are fairly robust with respect to small misspecifications of the cell shapes, they are attractive candidates for the data at hand. In Section \ref{simulation_sec}, we also pointed out that computing the persistence diagram via the tessellation vertices did not improve the testing power -- at least for the models considered in the simulation study. However, in other scenarios, the tests based on the centroids could be outperformed when replacing the cell centroids by other point clouds extracted from the slices.

Third, as pointed out in Section \ref{int_sec}, instead of computing test statistics related to topological characteristics of 2D slices, one may also attempt to reconstruct faithfully certain topological characteristics of the 3D data set. However, this construction relies on potentially delicate density and transversality conditions. It is interesting to investigate whether the longitudinal and cross-sectional tests from our framework also become more powerful under such conditions.

Finally, it would be exciting to generalize to higher-dimensional situations such as snapshots of 3D data evolving in time. The main challenge here lies in generalizing the concept of $M$-bounded features which currently relies on the duality in two dimensions.

\section*{Acknowledgements}
	We thank  J.~G.~L\'opez for providing the data. AC is supported by the Netherlands Organisation for Scientific Research (NWO) through grant 613.009.102. CH acknowledges the financial support of the CogniGron research center and the Ubbo Emmius Funds (University of Groningen). MV acknowledges the financial support of the Research Program of the Materials innovation institute (M2i) (\url{www.m2i.nl}) supported by the Dutch government.

\phantomsection
\addcontentsline{toc}{section}{References}
\bibliography{lit}

\begin{thebibliography}{30}
\expandafter\ifx\csname natexlab\endcsname\relax\def\natexlab#1{#1}\fi
\providecommand{\url}[1]{\texttt{#1}}
\providecommand{\href}[2]{#2}
\providecommand{\path}[1]{#1}
\providecommand{\DOIprefix}{doi:}
\providecommand{\ArXivprefix}{arXiv:}
\providecommand{\URLprefix}{URL: }
\providecommand{\Pubmedprefix}{pmid:}
\providecommand{\doi}[1]{\href{http://dx.doi.org/#1}{\path{#1}}}
\providecommand{\Pubmed}[1]{\href{pmid:#1}{\path{#1}}}
\providecommand{\bibinfo}[2]{#2}
\ifx\xfnm\relax \def\xfnm[#1]{\unskip,\space#1}\fi
\bibitem[{Amini et~al.(2013)Amini, Boissonnat and Memari}]{amini}
\bibinfo{author}{Amini, O.}, \bibinfo{author}{Boissonnat, J.D.},
  \bibinfo{author}{Memari, P.}, \bibinfo{year}{2013}.
\newblock \bibinfo{title}{Geometric tomography with topological guarantees}.
\newblock \bibinfo{journal}{Discrete Comput. Geom.} \bibinfo{volume}{50},
  \bibinfo{pages}{821--856}.
\bibitem[{Baryshnikov and Yukich(2005)}]{barysh}
\bibinfo{author}{Baryshnikov, Y.}, \bibinfo{author}{Yukich, J.E.},
  \bibinfo{year}{2005}.
\newblock \bibinfo{title}{Gaussian limits for random measures in geometric
  probability}.
\newblock \bibinfo{journal}{Ann. Appl. Probab.} \bibinfo{volume}{15},
  \bibinfo{pages}{213--253}.
\bibitem[{Bickel and Wichura(1971)}]{bickel}
\bibinfo{author}{Bickel, P.J.}, \bibinfo{author}{Wichura, M.J.},
  \bibinfo{year}{1971}.
\newblock \bibinfo{title}{Convergence criteria for multiparameter stochastic
  processes and some applications}.
\newblock \bibinfo{journal}{Ann. Math. Statist.} \bibinfo{volume}{42},
  \bibinfo{pages}{1656--1670}.
\bibitem[{Billingsley(1999)}]{billingsley}
\bibinfo{author}{Billingsley, P.}, \bibinfo{year}{1999}.
\newblock \bibinfo{title}{Convergence of Probability Measures}.
\newblock \bibinfo{edition}{Second} ed., \bibinfo{publisher}{John Wiley \&
  Sons}, \bibinfo{address}{New York}.
\bibitem[{Biscio et~al.(2020)Biscio, Chenavier, Hirsch and Svane}]{svane}
\bibinfo{author}{Biscio, C.A.N.}, \bibinfo{author}{Chenavier, N.},
  \bibinfo{author}{Hirsch, C.}, \bibinfo{author}{Svane, A.M.},
  \bibinfo{year}{2020}.
\newblock \bibinfo{title}{Testing goodness of fit for point processes via
  topological data analysis}.
\newblock \bibinfo{journal}{Electron. J. Stat.} \bibinfo{volume}{14},
  \bibinfo{pages}{1024--1074}.
\bibitem[{B{\l}aszczyszyn et~al.(2019)B{\l}aszczyszyn, Yogeshwaran and
  Yukich}]{yogeshCLT}
\bibinfo{author}{B{\l}aszczyszyn, B.}, \bibinfo{author}{Yogeshwaran, D.},
  \bibinfo{author}{Yukich, J.E.}, \bibinfo{year}{2019}.
\newblock \bibinfo{title}{Limit theory for geometric statistics of point
  processes having fast decay of correlations}.
\newblock \bibinfo{journal}{Ann. Probab.} \bibinfo{volume}{47},
  \bibinfo{pages}{835--895}.
\bibitem[{Chiu et~al.(2013)Chiu, Stoyan, Kendall and Mecke}]{daley2007II}
\bibinfo{author}{Chiu, S.N.}, \bibinfo{author}{Stoyan, D.},
  \bibinfo{author}{Kendall, W.S.}, \bibinfo{author}{Mecke, J.},
  \bibinfo{year}{2013}.
\newblock \bibinfo{title}{Stochastic geometry and its applications}.
\newblock Wiley Series in Probability and Statistics. \bibinfo{edition}{third}
  ed., \bibinfo{publisher}{John Wiley \& Sons, Ltd., Chichester}.
\bibitem[{Cohen-Steiner et~al.(2006)Cohen-Steiner, Edelsbrunner and
  Morozov}]{vine}
\bibinfo{author}{Cohen-Steiner, D.}, \bibinfo{author}{Edelsbrunner, H.},
  \bibinfo{author}{Morozov, D.}, \bibinfo{year}{2006}.
\newblock \bibinfo{title}{Vines and vineyards by updating persistence in linear
  time}, in: \bibinfo{booktitle}{Computational Geometry ({SCG}'06)}.
  \bibinfo{publisher}{ACM, New York}, pp. \bibinfo{pages}{119--126}.
\bibitem[{Davydov and Zitikis(2008)}]{davydov}
\bibinfo{author}{Davydov, Y.}, \bibinfo{author}{Zitikis, R.},
  \bibinfo{year}{2008}.
\newblock \bibinfo{title}{On weak convergence of random fields}.
\newblock \bibinfo{journal}{Ann. Inst. Statist. Math.} \bibinfo{volume}{60},
  \bibinfo{pages}{345--365}.
\bibitem[{Divol and Polonik(2019)}]{divol}
\bibinfo{author}{Divol, V.}, \bibinfo{author}{Polonik, W.},
  \bibinfo{year}{2019}.
\newblock \bibinfo{title}{On the choice of weight functions for linear
  representations of persistence diagrams}.
\newblock \bibinfo{journal}{J. Appl. Comput. Topol.} \bibinfo{volume}{3},
  \bibinfo{pages}{249--283}.
\bibitem[{Edelsbrunner and Harer(2010)}]{edHar}
\bibinfo{author}{Edelsbrunner, H.}, \bibinfo{author}{Harer, J.},
  \bibinfo{year}{2010}.
\newblock \bibinfo{title}{Computational Topology}.
\newblock \bibinfo{publisher}{American Mathematical Society, Providence, RI}.
\bibitem[{Eichelsbacher et~al.(2015)Eichelsbacher, Rai\v{c} and
  Schreiber}]{raic}
\bibinfo{author}{Eichelsbacher, P.}, \bibinfo{author}{Rai\v{c}, M.},
  \bibinfo{author}{Schreiber, T.}, \bibinfo{year}{2015}.
\newblock \bibinfo{title}{Moderate deviations for stabilizing functionals in
  geometric probability}.
\newblock \bibinfo{journal}{Ann. Inst. Henri Poincar\'e Probab. Stat.}
  \bibinfo{volume}{51}, \bibinfo{pages}{89--128}.
\bibitem[{Gal{\'a}n~L{\'o}pez and Kestens(2021)}]{galan}
\bibinfo{author}{Gal{\'a}n~L{\'o}pez, J.}, \bibinfo{author}{Kestens, L.A.},
  \bibinfo{year}{2021}.
\newblock \bibinfo{title}{A multivariate grain size and orientation
  distribution function: derivation from electron backscatter diffraction data
  and applications}.
\newblock \bibinfo{journal}{J. Appl. Crystallogr.} \bibinfo{volume}{54}.
\bibitem[{Hahn and Lorz(1993)}]{hahn}
\bibinfo{author}{Hahn, U.}, \bibinfo{author}{Lorz, U.}, \bibinfo{year}{1993}.
\newblock \bibinfo{title}{Stereological model tests for the spatial
  {P}oisson-{V}oronoi tessellation {II}}.
\newblock \bibinfo{journal}{Acta Stereol.} \bibinfo{volume}{12},
  \bibinfo{pages}{131--140}.
\bibitem[{Heinrich(2015)}]{ripley}
\bibinfo{author}{Heinrich, L.}, \bibinfo{year}{2015}.
\newblock \bibinfo{title}{Gaussian limits of empirical multiparameter
  {$K$}-functions of homogeneous {P}oisson processes and tests for complete
  spatial randomness}.
\newblock \bibinfo{journal}{Lith. Math. J.} \bibinfo{volume}{55},
  \bibinfo{pages}{72--90}.
\bibitem[{Krebs and Hirsch(2022)}]{krebs}
\bibinfo{author}{Krebs, J.T.N.}, \bibinfo{author}{Hirsch, C.},
  \bibinfo{year}{2022}.
\newblock \bibinfo{title}{Functional central limit theorems for persistent
  {B}etti numbers on cylindrical networks}.
\newblock \bibinfo{journal}{Scand. J. Stat.} \bibinfo{volume}{49},
  \bibinfo{pages}{forthcoming}.
\bibitem[{Last and Penrose(2016)}]{poisBook}
\bibinfo{author}{Last, G.}, \bibinfo{author}{Penrose, M.D.},
  \bibinfo{year}{2016}.
\newblock \bibinfo{title}{Lectures on the Poisson Process}.
\newblock \bibinfo{publisher}{Cambridge University Press},
  \bibinfo{address}{Cambridge}.
\bibitem[{Madej(2017)}]{madej}
\bibinfo{author}{Madej, L.}, \bibinfo{year}{2017}.
\newblock \bibinfo{title}{Digital/virtual microstructures in application to
  metals engineering--a review}.
\newblock \bibinfo{journal}{Arch. Civ. Mech. Eng.} \bibinfo{volume}{17},
  \bibinfo{pages}{839--854}.
\bibitem[{McLeish(1974)}]{mcleish}
\bibinfo{author}{McLeish, D.L.}, \bibinfo{year}{1974}.
\newblock \bibinfo{title}{Dependent central limit theorems and invariance
  principles}.
\newblock \bibinfo{journal}{Ann. Probab.} \bibinfo{volume}{2},
  \bibinfo{pages}{620--628}.
\bibitem[{Peccati and Taqqu(2011)}]{peccatitaqqu}
\bibinfo{author}{Peccati, G.}, \bibinfo{author}{Taqqu, M.S.},
  \bibinfo{year}{2011}.
\newblock \bibinfo{title}{Wiener Chaos: Moments, Cumulants and Diagrams}.
\newblock \bibinfo{publisher}{Springer, Milan; Bocconi University Press,
  Milan}.
\bibitem[{Penrose(2003)}]{penrose}
\bibinfo{author}{Penrose, M.D.}, \bibinfo{year}{2003}.
\newblock \bibinfo{title}{Random Geometric Graphs}.
\newblock \bibinfo{publisher}{Oxford University Press},
  \bibinfo{address}{Oxford}.
\bibitem[{Penrose and Yukich(2001)}]{yukCLT}
\bibinfo{author}{Penrose, M.D.}, \bibinfo{author}{Yukich, J.E.},
  \bibinfo{year}{2001}.
\newblock \bibinfo{title}{Central limit theorems for some graphs in
  computational geometry}.
\newblock \bibinfo{journal}{Ann. Appl. Probab.} \bibinfo{volume}{11},
  \bibinfo{pages}{1005--1041}.
\bibitem[{Pirgazi(2019)}]{pirgazi}
\bibinfo{author}{Pirgazi, H.}, \bibinfo{year}{2019}.
\newblock \bibinfo{title}{On the alignment of {3D EBSD} data collected by
  serial sectioning technique}.
\newblock \bibinfo{journal}{Mater. Charact.} \bibinfo{volume}{152},
  \bibinfo{pages}{223--229}.
\bibitem[{Reani and Bobrowski(2021)}]{reani}
\bibinfo{author}{Reani, Y.}, \bibinfo{author}{Bobrowski, O.},
  \bibinfo{year}{2021}.
\newblock \bibinfo{title}{Cycle registration in persistent homology with
  applications in topological bootstrap}.
\newblock \bibinfo{journal}{arXiv preprint arXiv:2101.00698} .
\bibitem[{Redenbach et~al.(2012)Redenbach, Shklyar and Andr\"a}]{stiff}
\bibinfo{author}{Redenbach, C.}, \bibinfo{author}{Shklyar, I.},
  \bibinfo{author}{Andr\"a, H.}, \bibinfo{year}{2012}.
\newblock \bibinfo{title}{Laguerre tessellations for elastic stiffness
  simulations of closed foams with strongly varying cell sizes}.
\newblock \bibinfo{journal}{Int. J. Eng. Sci.} \bibinfo{volume}{50},
  \bibinfo{pages}{70--78}.
\bibitem[{Salch et~al.(2021)Salch, Regalski, Abdallah, Suryadevara, Catanzaro
  and Diwadkar}]{fmri}
\bibinfo{author}{Salch, A.}, \bibinfo{author}{Regalski, A.},
  \bibinfo{author}{Abdallah, H.}, \bibinfo{author}{Suryadevara, R.},
  \bibinfo{author}{Catanzaro, M.J.}, \bibinfo{author}{Diwadkar, V.A.},
  \bibinfo{year}{2021}.
\newblock \bibinfo{title}{From mathematics to medicine: A practical primer on
  topological data analysis ({TDA}) and the development of related analytic
  tools for the functional discovery of latent structure in {fMRI} data}.
\newblock \bibinfo{journal}{PLOS ONE} \bibinfo{volume}{16},
  \bibinfo{pages}{1--33}.
\bibitem[{Tewari and Gokhale(2001)}]{tewari}
\bibinfo{author}{Tewari, A.}, \bibinfo{author}{Gokhale, A.M.},
  \bibinfo{year}{2001}.
\newblock \bibinfo{title}{Estimation of three-dimensional grain size
  distribution from microstructural serial sections}.
\newblock \bibinfo{journal}{Mat. Charact.} \bibinfo{volume}{46},
  \bibinfo{pages}{329--335}.
\bibitem[{Vittorietti et~al.(2020)Vittorietti, Kok, Sietsma, Li and
  Jongbloed}]{vittorietti}
\bibinfo{author}{Vittorietti, M.}, \bibinfo{author}{Kok, P.J.J.},
  \bibinfo{author}{Sietsma, J.}, \bibinfo{author}{Li, W.},
  \bibinfo{author}{Jongbloed, G.}, \bibinfo{year}{2020}.
\newblock \bibinfo{title}{General framework for testing {P}oisson-{V}oronoi
  assumption for real microstructures}.
\newblock \bibinfo{journal}{Appl. Stoch. Models Bus. Ind.}
  \bibinfo{volume}{36}, \bibinfo{pages}{604--627}.
\bibitem[{Wasserman(2018)}]{wasserman}
\bibinfo{author}{Wasserman, L.}, \bibinfo{year}{2018}.
\newblock \bibinfo{title}{Topological data analysis}.
\newblock \bibinfo{journal}{Annu. Rev. Stat. Appl.} \bibinfo{volume}{5},
  \bibinfo{pages}{501--535}.
\bibitem[{Yoo et~al.(2016)Yoo, Kim, Ahn and Ye}]{brain}
\bibinfo{author}{Yoo, J.}, \bibinfo{author}{Kim, E.Y.}, \bibinfo{author}{Ahn,
  Y.M.}, \bibinfo{author}{Ye, J.C.}, \bibinfo{year}{2016}.
\newblock \bibinfo{title}{Topological persistence vineyard for dynamic
  functional brain connectivity during resting and gaming stages}.
\newblock \bibinfo{journal}{J. Neurosci. Methods} \bibinfo{volume}{267},
  \bibinfo{pages}{1--13}.

\end{thebibliography}
\newpage

\setcounter{page}{1}
\section*{Supplementary material}

\setcounter{section}{0}
\renewcommand{\thesection}{\Alph{section}}
\section{Proof of Theorem \ref{clt_thm}}
\label{scal_sec}

%
%
To prove the scalar CLT stated in Theorem \ref{clt_thm}, note that the problem under consideration is very close to the general framework of Poisson-based test statistics as presented in \citet[Theorem 3.1]{yukCLT}. However, since in \eqref{eq:XX} we implement minus-sampling as an edge correction, the test statistic $T_n$ is not a pure Poisson functional, and therefore \citet[Theorem 3.1]{yukCLT} cannot be applied directly. However, the general methodology is sufficiently flexible to deal also with the present situation. In fact, the necessary modifications are essentially described in \cite{krebs}. Nevertheless, to make the manuscript self-contained, we reproduce here the most important steps. Without loss of generality, we may assume that the score function $\xi$ be bounded above by 1. Henceforth, $c_1,\,c_2,\,C\dots$ denote positive universal constants which may vary from line to line within the same equation.

%
%
Let $Q(z_1), \dots, Q(z_{\kn})$ denote the lattice cubes intersecting $\Qn$. Therefore, $k_n$ is of order $n^{p'}$. We assume henceforth that the $z_i \in \Z^{p'}$, $i \le \kn$ appear in the lexicographic order $\le_{\ms{lex}}$. Let 
$$\GG_i := \GG_{z_i} := \s\big(\PP \cap \bigcup\nolimits_{z \le_{\ms{lex}} z_i}Q(z)\big)$$
denote the $\s$-algebra generated by the Poisson points in cubes $Q(z)$ with $z \le_{\ms{lex}} z_i$.
Now, we put 
$$\mc D_{i, n} := \E[T_n\ba \GG_{i}] - \E[T_n\ba \GG_{{i - 1}}] ,$$
so that $T_n$ admits the martingale-difference decomposition
$T_n - \E[T_n] = \sum_{i \le \kn}\mc D_{i, n}.$

The key ingredient for the proof of Theorem \ref{clt_thm} is the martingale CLT~\cite[Theorem 2.3]{mcleish}. We restate it here following~\citet[Theorem 2.10]{penrose}, when specialized to mean-zero martingales.
\bet\label{mcl_thm}
			Suppose that $\{M_{i, n}\}_{1 \le i \le n}$ is a mean-zero martingale for each $n \ge 1$ and set $Y_{i, n} := M_{i, n} - M_{i - 1, n}$ with $M_{0, n} = 0$. Suppose that
\begin{enumerate}[label=(\alph*),ref=(\alph*)]
	\item~\label{a} $\sup_{n\ge 1} \E[ \max_{i \le n} Y_{i, n} ^2 ] < \infty $;
	\item~\label{b} $  \max_{i \le n} |Y_{i, n}| \to 0 $ in probability as $n\to\infty$;
						\item~\label{c} $ \sum_{i \le n } Y_{i, n}^2 \to \wt \s^2 $ in $L^1$ for some $\wt\s^2 > 0$.
							\end{enumerate}
			        Then, $M_{n, n} \to \NN(0, \wt\s^2)$ in distribution.
				\end{theorem}
Henceforth, we verify conditions~\ref{a}--\ref{c} for $Y_{i, \kn} := \kn^{-1/2}\mc D_{i, n}$. To that end, we need a stabilization result and a moment bound. To state them precisely, we provide a more conceptual description of $\mc D_{i, n}$. 
Indeed, let $\PP_n'$ be an independent copy of $\PP_n$. Then, for $z \in \Z^p$ let $T_{z, n}$ denote the test statistic $T_n$ computed on the basis of
$$\PP_{z, n} :=\big(\PP_n\sm Q(z)\big) \cup \big(\PP_n'\cap  Q(z)\big)$$
instead of $\PP_n$.  Furthermore, set $\De_{z, n} := T_n - T_{z, n}$. Then, $\E[T_n\ba \GG_{{i - 1}}] = \E[T_{z_i, n} \ba \GG_{i}]$  so that 
$\mc D_{i, n} = \E[\De_{z_i,n }\ba \GG_{i} ]$. We will use the following Lemma to prove~\ref{a}--\ref{c}. Its proof can be found on page~\pageref{pg:lemma}.
%
%
\begin{lemma}[Uniformly bounded moments]\label{unif_lem}
	        It holds that $\sup_{n \ge 1} \sup_{i \le \kn}\E\big[\De_{z_i, n}^4\big] < \infty$. Thus,
		\begin{align}\label{E:UniformBounded1}
			        \sup_{n \ge 0} \sup_{i \le \kn}\E\big[\mc D_{i, n}^4\big] < \infty.
				 \end{align}
\end{lemma}
To verify condition~\ref{c} we are going to show
%
%
\bel[Stabilization]
\label{stab_lem}
	Let $z \in \Z^{p'}$. Then, $\De_{z, n} = \De_{z, n'}$ for every $n, n' \ge 1$ such that $Q_{4R(z, n)} \su Q_n(z)$ and $Q_{4R(z, n')}(z) \su Q_{n'}(z)$. We write $\De_{z, \ff}$  for the common value.
\enl
The proof of this Lemma is deferred to page~\pageref{pg:stab_lemma}.
\begin{proof}[Proof of Theorem~\ref{mcl_thm}] 
We are going to verify conditions~\ref{a}--\ref{c}.
For condition~\ref{a}:
\begin{align*}
	&\sup_{n\ge 1} \E[ \max_{i \le \kn} Y_{i, n} ^2 ] = \sup_{n\ge 1} \ \kn^{-1} \ \E[ \max_{i \le \kn}  \mc D_{i, n} ^2 ] \le \sup_{n\ge 1} \ \kn^{-1} \ \sum_{i\le \kn} \E\big[ \mc D_{i, n}^2 \big] = O(1).
\end{align*}
In the last step we used~\eqref{E:UniformBounded1}. 
	To verify condition~\ref{b}, note that 
\begin{align*}
 \quad  \E[ \max_{i \le \kn} Y_{i, n} ^4 ] = \kn^{-2} \ \E[ \max_{i \le \kn}  \mc D_{i, n}^4 ] \le \kn^{-2} \sum_{i\le \kn} \E\big[ \mc D_{i, n}^4 \big] = O(k_n^{-1})
\end{align*}
	using again \eqref{E:UniformBounded1}.
Hence, $\max_{i \le \kn} Y_{i, n}$ converges to $0$ in $L^4$ and therefore also in probability.
 We elucidate how Lemma \ref{stab_lem} allows us to verify condition~\ref{c} of Theorem \ref{mcl_thm}.

First, we deduce from Lemma \ref{unif_lem} that the sum 
$\kn^{-1}\sum_{i\co Q(z_i) \not \su Q_{n - 2\sqrt n}} \mc D_{i, n}^2$
tends to 0 in $L^1$. Next, since the sequence of random variables $\bi\{\E\bi[ \De_{z, \ff}^4  \ba \GG_z\bi]\bi\}_{z \in \Z^{p'}}$ is stationary, we deduce from the Cauchy-Schwarz inequality that
$$\kn^{-1}\E\Big[\sum_{i\le \kn} \E\bi[ \De_{z_i, \ff} \one_{\{R(z_i) > \sqrt n/4\}}\bi \ba \GG_i\bi]^2\Big] \le  \sqrt{\E\bi[ \De_{o, \ff}^4\bi]} \sup_{i \le \kn}\sqrt{\P(R(z_i) > \sqrt n/4)}.$$
	{By exponential stabilization and Lemma~\ref{unif_lem}}, the right-hand side tends to $0$ as $n \to \ff$. Similarly, we can also bound
$$\kn^{-1}\E\Big[\sum_{i\le \kn} \E\bi[ \De_{z_i, n} \one_{\{R(z_i) > \sqrt n/4\}}\bi \ba \GG_i\bi]^2\Big] \le \kn^{-1}\sum_{i\le \kn} \E\bi[ \De_{z_i, n}^2 \one_{\{R(z_i) > \sqrt n/4\}}\bi].$$
Moreover, Lemma \ref{stab_lem} yields that almost surely,
	$$\sum_{i\co Q(z_i) \su Q_{n - 2\sqrt n}} \kn^{-1}\E\bi[ \De_{z_i, n} \one_{\{R(z_i) \le \sqrt n/4\}} \ba \GG_i\bi]^2 = \sum_{i\co Q(z_i) \su Q_{n - 2\sqrt n}} \kn^{-1}\E\bi[ \De_{z_i, \ff} \one_{\{R(z_i) \le \sqrt n/4\}} \ba \GG_i\bi]^2.$$
Hence, $\kn^{-1}\sum_{i\le \kn} \E\bi[ \De_{z_i, n}  \ba \GG_i\bi]^2$ and $\kn^{-1}\sum_{i\le \kn} \E\bi[ \De_{z_i, \ff}  \ba \GG_i\bi]^2$
	agree up to an $L^1$-negligible error. To conclude the verification of condition~\ref{c}, we observe that $\bi\{\E\bi[ \De_{z, \ff}  \ba \GG_z\bi]^2\bi\}_{z \in \Z^{p'}}$ is stationary, so that we may apply the ergodic theorem as in the proof of \cite[Theorem 3.1]{yukCLT}.
\end{proof}

It remains to prove Lemmas \ref{unif_lem} and \ref{stab_lem}. The key observation for the proof of Lemma \ref{unif_lem} is that exponential stabilization allows to pass from $\PP_n$ to $\PP_{z_i, n}$ affecting only very few trajectories.

\bep[Proof of Lemma \ref{unif_lem}]~\label{pg:lemma}
%
%
To simplify notation, we write $\De_{i, n}$ instead of $\De_{z_i, n}$. Since  $\xi$ is bounded above by 1, we deduce that $|\De_{i, n}|$ is at most  the number of points in the $M$-bounded persistence diagrams on $\PP_n$ and $\PP_{z_i, n}$. First, the number of $M$-bounded clusters is bounded above by the number of points in the underlying point cloud so that $|\De_{i, n}| \le |\PP_n| + |\PP_{i, n}|$ in degree 0. In degree 1, we note when a new hole is born at level $r$, then there are two points of the point cloud at distance $2r$, and conversely each such pair gives rise to at most two new $M$-bounded holes. Thus, 
\begin{align}
	\label{deb_eq}
	|\De_{i, n}| \le 2|\PP_n|^2 + 2|\PP_{i, n}|^2.
\end{align}
Since $\PP_n$ and $\PP_{i, n}$ have the same distribution, we deduce from the Cauchy-Schwarz inequality that
\begin{equation}\label{eq:CS-exp-radius}\E[|\De_{i, n}|^k \one_{\{R(z_i) > \sqrt n /4\}}] \le 2^k\sqrt{\E[|\PP_n|^{4k}]}\sqrt{\P(R(z_i)  \ge \sqrt n/4 )}.\end{equation}
By exponential stabilization,  right-hand side {is negligible} as $n \to \ff$.  Hence, we may concentrate on the event that $R(z_i, n)\le \sqrt n /4$. In particular, if additionally, $Q(z_i) \not \su Q_{n - \sqrt n /2}$, then condition~\ref{S1} implies that $\De_{i, n} = 0$. It therefore remains to treat the case $Q(z_i) \su Q_{n - \sqrt n /2}$.

%
%
Now, the definition of the stabilization radius $R({z_i})$ gives that
\begin{align}
	\label{dnxz_eq}
	|\dnz| \le \dexx  + \De_{i, z_i, n}' ,
\end{align}
 where $\dexx$ denotes the total number of $M$-bounded $q$-features in $\XX_n$ associated with Poisson points in $Q_{2R({z_i}, n)}(z_i)$.   The quantity $\De_{i, z_i, n}'$ is defined in the same manner, except that we replace $\XX_n$ by 
$$
\XX_{z_i, n} :=  \{\TT(P_i, \PP_{{z_i},n})\}_{P_i \in \PP_{z_i,\,n - \sqrt n}}.
 $$ Indeed, by condition~\ref{S1} changing $\PP_n$ in $Q(z_i)$ can only change features that are associated with a trajectory generated by $P \in \PP \cap Q_{R(z_i)}(z_i)$. By condition~\ref{S2} and the $M$-boundedness, any such feature consists entirely of trajectories associated with points in $P \in \PP \cap Q_{2R(z_i)}(z_i)$. 

Since $\dexx$ and $\De_{i, z_i, n}'$ have the same distribution, it suffices to bound the moments of $\dexx$. To achieve this goal, we first claim that   $|\dexx| \le 2\PP_n( Q_{2R(z_i, n)}(z_i))^2$. We present the reasoning for the $M$-bounded holes, noting that the case of $M$-bounded clusters is similar but easier. 
Now, as in the derivation of \eqref{deb_eq}, we note that any pair of points from $\PP \cap Q_{2R(z_i)}(z_i)$ can give rise to at most two $M$-bounded holes, which therefore leads to the bound $\dexx \le 2\PP(Q_{2R(z_i, n)}(z_i))^2$.  Thus, by the Cauchy-Schwarz inequality,
\begin{align*}
	\E[(\dexx)^k] \le 	2^k\E\big[\PP( Q_{2R(z_i, n)}(z_i))^{2k}\big] 
	&\le 2^k\sum_{m \ge 1} \sqrt{\E\big[\PP(  Q_{2m}(z_i))^{4k}\big]} \sqrt{\P( R(z_i, n) = m)}\\
	&\le 2^k\sqrt{\E\big[\PP(Q)^{4k}\big]} \sum_{m \ge 1}(2m)^{4kd + d}\sqrt{\P( R(z_i, n) = m)}, 
\end{align*}
and the stretched exponential decay of $R(z_i, n)$ guarantees the finiteness of the right-hand side.
\enp

\bep[Proof of Lemma \ref{stab_lem}]\label{pg:stab_lemma}
We first  apply property {\bf (S1)} with $\AA = \PP_n \sm Q_{R(z,n)}(z)$ and $\BB = \PP_n \cap Q(z)$. Thus, all trajectories in the symmetric difference $\XX_\De$ between $\XX_n$ and $\XX_{z, n}$ are associated with points $P \in \PP_n \cap Q_{R(z, n)}(z)$. Second, we apply property {\bf (S2)} with $\AA = \PP_n \sm Q_{2R(z, n)}(z)$ and $\BB = \PP \cap Q$. Thus, any $M$-bounded feature involving a trajectory from $\XX_\De$ is associated with a point $P \in \PP_n \cap Q_{2R(z, n)}(z)$. Finally,  we apply property {\bf (S3)} with $\AA = \PP_n \sm Q_{4R(z, n)}(z)$. Thus, such trajectories remain unchanged, once $n$ is so large that $Q_{4R(z, n)}(z) \su \Qn$.
\enp

\section{Proof of Theorem \ref{fclt_thm}}
\label{tight_sec}
The proof of Theorem \ref{fclt_thm} combines arguments appearing in \cite{svane, krebs}. On the one hand, similarly to \cite{krebs} our null hypothesis is based on a Poisson point process, so that we can build on the martingale approach from \cite{yukCLT}. Moreover, we heavily rely on stabilization properties of the trajectory construction rule, which are reminiscent of the stabilizing networks encountered in \cite{krebs}.

On the other hand, the sampling window $Q_n$ now grows not just in one dimension so that percolation effects need to be controlled. Hence, we restrict to $M$-bounded features as in \cite{svane}. Moreover, we face the additional challenge that the point cloud giving rise to the persistence diagram arises from the underlying Poisson process through a possibly complicated construction rule. In order to establish the continuity properties appearing in the standard tightness criteria, the boundedness of the factorial moment densities from condition \eqref{mom_eq} is crucial.

Henceforth, we explain the proof for the feature dimension $q = 1$, i.e., for $M$-bounded holes. The case $q = 0$ is in large parts parallel but less complicated since here the persistent Betti numbers are not a bivariate but a univariate process. Moreover, we always tacitly assume that features are $M$-bounded for some fixed $M > 0$, and hence will suppress the dependence on $M$ to simplify notation.
Proving  asymptotic normality on a functional level consists of two steps: 1) multi-variate normality of the marginals and 2) tightness. First, invoking the Cram\'er-Wold device~\cite[Theorem 7.7]{billingsley}, the multivariate normality becomes a consequence of Theorem \ref{clt_thm}. 

%
%
In order to prove tightness, we rely on the Chentsov-type tightness condition from \cite{bickel}. More precisely, for a block 
$$E :=  \Eb \ti \Ed :=  [b_-, b_+] \ti [d_-, d_+] \su  [0, \t]^2 ,$$ 
we write 
$$\b_n(E) := \PD_n(E) = \f1H\sum_{h \in \Xi_H}\#\big\{i \ge 1\co (B_i(h), D_i(h)) \in \Eb \ti \Ed \big\} $$
for the contribution of the features in $E$ and set $\bar\b_n(E) := \b_n(E) - \E[\b_n(E)]$. Then, we need to show that 
\begin{align}
\label{chentsov_eq}
n^{-2p'}\E[\bar\b_n(E)^4] \le C |E|^{1 + \e}.
\end{align}

%
%
To achieve this goal, we proceed similarly to the proof of \cite[Theorem 2]{krebs}. First, in Proposition \ref{grid_sec} proven in Subsection~\ref{grid_sec}, we argue that it suffices to verify \eqref{chentsov_eq} for \emph{$n$-big blocks}, i.e., blocks $E =  \Eb \ti \Ed$ satisfying $|E| \ge n^{-2p'}$. 

%
\bepr[Reduction to grid]
\label{grid_prop}
The processes $\{\bnn\}_{b, d}$ are tight in the Skorokhod topology if condition \eqref{chentsov_eq} is satisfied for all $n \ge 1$ and all $n$-big $E \su [0, \t]^2$.
\enpr

Second, in Proposition \ref{md_sec} proven in Subsection~\ref{md_sec}, we leverage a martingale-difference argument in order to derive the following key variance and cumulant bounds in the vein of~\cite[Proposition 8]{krebs}. Henceforth, $c(X_1,\,\dots,\,X_d)$ denotes the joint cumulant of the random variables $(X_1,\,\dots,\,X_d)$, while for a single random variable $X$ the quantity $c^4(X) = c(X, X, X, X)$ denotes the fourth cumulant of $X$.

%
%
\bepr[Variance and cumulant bound -- large blocks]
\label{var_prop}
It holds that 
$$\sup_{n \ge 1}\sup_{\substack{E \su [0, \t^2] \\ E \text{ is $n$-big}}}\f{\Var\big(\b_n(E)\big) +  c^4\big(\b_n(E)\big)}{n^{p'} |E|^{5/8}} < \ff.$$
\enpr

%
%
We are now ready to prove Theorem~\ref{fclt_thm}.
\bep[Proof of Theorem~\ref{fclt_thm}]
Decomposing into a variance and a cumulant contribution gives that
$$n^{-2p'}\E[\bar\b_n(E)^4] = 3n^{-2p'}\Var(\b_n^2(E))^2 + n^{-2p'} c^4(\b_n(E)).$$
Hence, applying Proposition \ref{var_prop} shows that the right-hand side is at most $3 C|E|^{5/4} + C|E|^{1/2 + 5/8}$ for a suitable $C > 0$, thereby concluding the proof.
\enp

\subsection{Reduction to grid}
\label{grid_sec} 
Proceeding in a similar vein as in \cite[Proposition 5]{krebs}, we now invoke a trick from \cite{davydov} allowing us to restrict to blocks with end points in a sufficiently fine grid of the form $G_n := n^{-p'}\Z^2 $, $n \ge 1$. More precisely, by \cite[Theorem 16.8]{billingsley} it remains to bound the continuity modulus
$$\omega'_\de(\bar\b_n) := \inf_\Gamma \max_{L\in\Gamma} \sup_{( b, d), ( b', d')\in L} |\bar\b_n^{  b, d} - \bar\b_n^{  b', d'}|, $$
where the infimum extends over all $\de$-grids $\Gamma$ in $ [0, \t]^2$.

%
%
For a $q$-simplex $\s = \{x_0, \dots, x_q\}$, we let 
$$r(\s) := \inf\bi\{t > 0\co \bigcap_{i \le q} B(x_i; t) \ne \es\bi\}$$
denote the \emph{filtration time} in the \v Cech filtration.
In order to prove Proposition \ref{grid_prop}, we need a refined upper bound on the probability that the \v Cech filtration time of a simplex $\s$ lies in a small interval. To that end, we rely on the following bound from the proof of \cite[Lemma 6.10]{divol}.

\bel[\v Cech filtration times are Lipschitz]
\label{divol_lem}
Let $ q \le p'$  and $X_0,\dots, X_q $ be iid uniform in a box $Q \su \R^{p'}$. Then, the function $r \mapsto \P(r(X_0, \dots, X_q) \le r)$ is Lipschitz.
\enl
%
%
\bep[Proof of Proposition \ref{grid_prop}]
To prove Proposition \ref{grid_prop} we need to show that for every $\e, \de > 0$ there exists $n_0 = n_0(\e, \de)$ such that almost surely
	        $$\sup_{n \ge n_0}n^{-p'/2}|\omega'_\de(\bar\b_n)- \omega'_\de(\bar\b_n|_{G_n})| \le \e.$$

Since the process $\b^{r, s}$ is increasing in $ r$ and $(-s)$, it suffices by \cite[Corollary 2]{davydov} to control
$$\E\big[ \b^{  r_2, s}_n - \b^{  r_1, s}_n \big] \quad\text{ and }\quad   \E\big[ \b^{  r, s_1}_n - \b^{  r, s_2}_n\big]
	$$
	for $  r_2 - r_1, s_2 - s_1 \in [0, n^{ -p'}]$.  We only tackle the second expression, as the arguments for the first are similar. 
	Now, let 
	$$E_n := \big\{\max_{i \le \kn}R(z_i, n) \le \sqrt n /4\big\}$$
	be the event that the stabilization radii at $z_i$, $i \le \kn$, are all at most $\sqrt n/4$. Furthermore, as in the derivation of \eqref{deb_eq}, we conclude that $\b^{  r, s_i}_n \le2 \PP( Q_n)^2$. Thus, by the Cauchy-Schwarz inequality,
$$\E\bi[(\b^{  r, s_1}_n - \b^{  r, s_2}_n)\one_{\{E_n^c\}}\bi] \le 2n\bi(\E\bi[\PP(Q_n)^4\bi]\bi)^{1/2} \P(E_n^c)^{1/2},$$
so that the stretched exponential decay of the stabilization radii shows that the right-hand side {is negligible for $n\to\infty$ similarly to~\eqref{eq:CS-exp-radius}}. Thus, we may henceforth work under the event $E_n$. 

Then, the increment $\b^{  r, s_1}_n - \b^{  r, s_2}_n$ is bounded by $\max_{h \in \Xi_H}g(h)$, where $g(h)$ is the number of $M$-bounded holes in $\XX_n(h)$ with death time in the interval $[s_1, s_2]$. Associating each feature with the 2-simplex causing causes its death, this quantity is bounded above by the number of 2-simplices with filtration time in the interval $[s_1, s_2]$.
Now, under the event $E_n$, we deduce from property {\bf (S3')} that all these simplices are contained in $Q_n$. Thus,  we invoke condition \eqref{mom_eq} to deduce that
\begin{align*}
	&\E\big[ (\b^{ r, s_1}_n - \b^{ r, s_2}_n)\one_{\{E_n \}}\big] \\
	&\quad\le \max_{h \in \Xi_H} \E\big[\# \{ \text{$2$-simplices } \s\in \PP_n(h)\co r(\s) \in [s_1, s_2]  \}\big] \\
	&\quad=\max_{h \in \Xi_H}\int_{Q_n'}\int_{Q_{2\t}'(x_0)^2} \one_{\big\{r(x_0, x_1, x_2) \in [s_1, s_2]\big\} }\rh(x_0, x_1, x_2) \d(x_1, x_2)\d x_0\\
	&\quad\le C_{\r, 3}\int_{Q_n'}\int_{Q_{2\t}'(x_0)^q} \one_{\big\{r(x_0, x_1, x_2) \in [s_1, s_2]\big\} } \d(x_1,x_2)\d x_0.
\end{align*}
Hence, an application of Lemma \ref{divol_lem} concludes the proof.
\enp

\subsection{Variance and cumulant bounds}
\label{md_sec}

%
%
As in the proof of Theorem \ref{clt_thm},
to establish the variance and cumulant bounds from Proposition \ref{var_prop}, we rely on the martingale-difference decompositions from \cite{yukCLT}. Now, we put 
$$\mc D_{i, n}(E) := \E[\b_n(E)\ba \GG_i ] - \E[\b_n(E)\ba \GG_{i - 1}] ,$$
so that $\bar\b_n(E)$ admits the decomposition 
$\bar\b_n(E) = \sum_{i \le \kn}\mc D_{i, n}(E).$

%
%
This decomposition already indicates how we should proceed for the variance computation. Indeed, 
$$\Var(\b_n(E)) = \sum_{i \le \kn} \Var(\mc D_{i, n}(E)),$$
so that as in the arguments in Lemma \ref{unif_lem}, we aim to establish a uniform control over the second moments of $\mc D_{i, n}(E)$ for $i \le \kn$ and $n \ge 1$ . However, now the bound also needs to reflect the block size $|E|$.
%
%
\bel[Moment bound]
	\label{diffBoundLem}
	It holds that
			$$\sup_{E \su [0, \t]^2}\sup_{n \ge 1}\max_{i \le \kn}   \f{\E[|\mc D_{i, n}(E)|^k]}{|E|^{11/16}} < \ff.$$
\enl

%
%
When passing to cumulants, we also need to control correlations. To that end, we define the distance 
\[
 \ms{dist}(X,\,Y):=\inf\{|x_i-y_j|:\, i\in I,\,j\in J\}.
\]
between points $X=\{x_i:\,i\in I\}\su \R^{p'}$, $Y=\{y_j:\,j\in J\}\su \R^{p'}$.

\begin{lemma}[Covariance bound]
	\label{covBoundLem}
	For every $r_1, r_2 \ge 1$ there exist $C_{r_1, r_2} C_{r_1, r_2}' > 0$ such that the following statement holds. Let $n \ge 1$ and $I_1, I_2 \su \{1, \dots, \kn\}$  with $|I_1| = r_1$, $I_2 = r_2$, and set          $X_1 = \prod_{i \in I_1}\mc D_{i, n}(E)$ and $X_2 = \prod_{j \in I_2 }\mc D_{j, n}(E)$. Then,
	$$\Cov\big(X_1, X_2\big) \le C_{r_1, r_2} \exp\big(-\dist(\{z_i\}_{i \in I_1},\{z_j\}_{j \in I_2})^{C_{r_1, r_2}'}\big)\sqrt{\E[X_1^4]\E[X_2^4]}.$$
\end{lemma}

Recall that the $z_i$'s denote the centers of the lattice cubes $Q_i$.
To ease notation, we henceforth write $\mc D_i$ for $\mc D_{i, n}(E)$ and $\b_n$ for $\b_n(E)$. The proof of Proposition \ref{var_prop} proceeds along the lines of \cite[Proposition 3]{krebs}. Nevertheless, to make the presentation self-contained, we include some details.

\bep[Proposition \ref{var_prop}]
The multilinearity of cumulants yields that
\begin{align}\label{Cum_LargeBlocks}
	c^4(\bar\b_n) \le \sum_{ i ,j ,k ,\ell\le \kn }a_{i, j, k, \ell} \ c\big(\mc D_i, \mc D_j, \mc D_k, \mc D_\ell\big),
\end{align}
where the $a_{i, j, k, \ell} \ge 1$ are suitable combinatorial coefficients such that these coefficients depend only on which of the indices $i, j, k, \ell$ are equal. In the following, up to a permutation of $(i, j, k, \ell)$  we need to distinguish three cases in bounding the right-hand side of \eqref{Cum_LargeBlocks}, where we set $\de := 1/64$: 
\begin{enumerate}[label=\Roman*.,ref=\Roman*.,leftmargin=*]
\im\label{cum_a} $\diam(\{z_i, z_j, z_k, z_\ell\}) < |E|^{-\de}$,
\im\label{cum_b} $\dist(z_i, \{z_j, z_k, z_\ell\}) \ge |E|^{-\de}$, 
\im\label{cum_c} $\dist(\{z_i, z_j\}, \{z_k, z_\ell\}) \ge |E|^{-\de}$ and $|z_i - z_j| \vee |z_k - z_\ell|\le |E|^{-\de}$.
\end{enumerate}
\begin{enumerate}[leftmargin=*]
\item[\ref{cum_a}] We need to derive an upper bound for
$$\sum_{i,j,k,\ell \co \ref{cum_a}}a_{i, j, k, \ell} \ c(\mc D_i, \mc D_j, \mc D_k, \mc D_\ell),$$
 where the subscript means that the summation extends over all indices satisfying condition \ref{cum_a}
We achieve this goal by invoking the H\"older inequality in order to exert control on individual cumulants. Indeed, together with the representation in \cite[Eq. 3.9]{raic}, we obtain that
\begin{align}
	\label{sing_cum_bound}
	\big|c(\mc D_i, \mc D_j, \mc D_k, \mc D_\ell)\big| &\le \sum_{L_1, \dots, L_r}a_{\{L_1, \dots, L_r\}}'
	\prod_{m \in L_1} \E[ |\mc D_m|^{|L_1|} ]^{1/|L_1|} \ \cdots \prod_{m \in L_r} \E[ |\mc D_m|^{|L_r|} ]^{1/|L_r|},
\end{align}
where  $L_1, \,\dots, L_r$ partitions $\{i, \dots, \ell\}$ and $a_{\{L_1, \dots, L_r\}}'$ are coefficients. Again, these coefficients  depend only on the type of the partition but not on the precise values of $i,j, k, \ell$.
Then, \eqref{sing_cum_bound} is at most   
\begin{align}
	\label{mmax_bound_eq}
	c \sup_{n \ge 1}\max_{\substack{k \le 4 \\ m \le \kn}}\E[|\mc D_m|^k]
\end{align}
for some $c > 0$ depending only on the chosen trajectory model. By Lemma \ref{diffBoundLem}, this expression is at most
$c'|E|^{11/16},$
where $c' > 0$ again depends only  on the underlying trajectory model. Hence,
$$\sum_{L_1,\, \dots,\, L_r}a_{i, j, k, \ell}  c(\mc D_i, \mc D_j, \mc D_k, \mc D_\ell)\le 2c'n|E|^{11/16}|E|^{-3p'\de},$$
and the right-hand side is in $O(n|E|^{5/8})$ for $|E|$ small.
 \item[\ref{cum_b}, \ref{cum_c}] Since the two cases are very similar, we provide only  the details for \ref{cum_b} Here, we want to bound the expression
	 $$\sum_{i,j,k,\ell \co \ref{cum_b}}a_{i, j, k, \ell} c(\mc D_i, \mc D_j, \mc D_k, \mc D_\ell).$$
To this end, we will use the idea presented in~\cite[Lemma 5.1]{barysh} for point measures to prove that
\begin{align}\label{eq:fourth_cum}
	c(\mc D_i, \mc D_j, \mc D_k, \mc D_\ell)  =\sum_{L_1, \dots, L_r} a_{\{L_1, \dots, L_r\}}'\Cov\Big(\mc D_i, \prod_{s \in L_1}\mc D_s\Big)\E\Big[\prod_{s \in L_2}\mc D_s\Big]\cdots \E\Big[\prod_{s \in L_r}\mc D_s\Big]
\end{align}
for some coefficients $a_{\{L_1, \dots, L_r\}}'$. Again, these coefficients depend only  on the type of the partition, but not on the precise values of $j, k, \ell$. 
		The proof of~\eqref{eq:fourth_cum} will be deferred to Section \ref{aux_sec}, since it follows closely that of~\cite{barysh}.
Hence, combining the moment bounds from Lemma \ref{diffBoundLem} with the covariance bounds from Lemma \ref{covBoundLem} concludes the proof.\qedhere
\end{enumerate}
\end{proof}
To prove the moment bound in Lemma \ref{diffBoundLem}, we can to a large extent follow the proof of Lemma \ref{unif_lem}. In order to avoid redundancies, we focus on the new steps needed to produce the block size $|E|$.

%
%
\bep[Proof of Lemma \ref{diffBoundLem}]
As in \eqref{dnxz_eq}, the definition of the stabilization radius $R(z_i, n)$ gives that
\begin{align}
	\label{dnx_eq}
	|\dnz| \le \max_{h \in \Xi_H} \dex + \max_{h \in \Xi_H}\De_{i, z_i, n}'(h) ,
\end{align}
 where $\dex$ denotes the total number of $M$-bounded $q$-features in $\XX_n(h)$ with birth and death time contained in $\Eb$ and $\Ed$ and that are contained in $Q_{2R({z_i, n})}(z_i)$. The quantity $\De_{i, z_i, n}$ is defined in the same manner, except that as before, we replace $\XX_n$ by $\XX_{z_i, n}$.

Thus, arguing as in the proof of Lemma \ref{unif_lem}, by the H\"older inequality,
\begin{align*}
	\E\Big[ \dex^k \Big] 
	&\le \E[\PP_n(Q_{2R(z_i, n)}(z_i))^{24k}]^{\f1{12}} \Big(\P\big(\dex \ge 1 \big)\Big)^{\f{11}{12}}\\
	&\le \E[\PP_n(Q_{R(z_i, n)}({z_i}))^{24k}]^{\f1{12}} |\Eh|^{\f{11}{12}}\P\big(\dex \ge 1 \big)^{\f{11}{12}}.
\end{align*}
Hence, it suffices to bound the probability $\P\big(\dex \ge 1 \big)$ uniformly over all $h, i$ and $n$. To achieve this goal, we first condition on the value of $R(z_i, n)$ to obtain that 
\begin{align*}
	\P\big(\dex \ge 1 \big) \le \sum_{m \ge 1} \P(m - 1 \le R(z_i, n) \le m)\P(\De_{z_i, 2m, n}''(h) \ge 1),
\end{align*}
where $\De_{z_i, 2m, n}''(h)$ denotes the total number of $M$-bounded $q$-features in $\XX_n(h)$ with birth and death time contained in $\Eb$ and $\Ed$, and that consist of trajectories associated with Poisson points in $Q_{2m}(z_i)$. Taking into account the stretched exponential decay of $R(z_i, n)$, it suffices to show that 
\begin{align}
	\label{d15_eq}
	\P(\De_{z_i, m, n}''(h) \ge 1) \le c m^{p'}|\Eb|^{3/4}|\Ed|^{3/4}
\end{align}
for some $c > 0$ depending only on the considered model. 

To establish the bound \eqref{d15_eq}, we first note that if $\De_{z_i, m, n}''(h) \ge 1$, then each $M$-bounded $q$-feature yields a pair of simplices $\s = \big(X_{i_0}(h), X_{i_1}(h)\big)$ and $\s' = \big(X_{i_0'}(h), X_{i_1'}(h), X_{i_2'}(h)\big)$ contained in $Q_{4m}(z_i)$ such that $r(\s) \in \Eb$ and $r(\s')\in \Ed$. {Essentially, $\sigma$ and $\sigma'$ are the simplices determining the birth and death of that particular feature}. Now, if $Q(\bar x_1), \dots, Q(\bar x_{m^{p'}})$ is a partition of $Q_m (z_i)$ into unit boxes, then we conclude that additionally that $\s \su Q_{2\t + 1}(\bar x_j)$ and $\s' \su Q_{2\t + 1}(\bar x_{j '})$ for some $j, j' \le m ^{p'}$. Thus, writing $\Ebd = \Eb \ti \Ed$ and using Markov's inequality,
\begin{align*}
	&\P(\De_{x_i, m}'(h)' \ge 1) \\
	&\le \sum_{j, j' \le m^{p'}}\P\big( (r(\s), r(\s'))\in \Ebd, r(\s') > r(\s) \text{ for some }\s \su Q_{2\t + 1}(\bar x_j), \s' \su Q_{2\t + 1}(\bar x_{j '})\big)\\
	&\le m^{2p'} \sup_{j, j' \le m^{p'}}\P\big( (r(\s), r(\s'))\in \Ebd, r(\s') > r(\s) \text{ for some }\s \su Q_{2\t + 1}(\bar x_j), \s' \su Q_{2\t + 1}(\bar x_{j '})\big)
\end{align*}
Finally, to bound the probability on the right-hand side, we need to distinguish on the number of vertices that $\s$ and $\s'$ have in common. We elaborate in detail on how to proceed if $\s \su \s'$, noting that the arguments in the other cases are very similar. Hence, $\s = (\{X_0, X_1\})$ and $\s' = (\{X_k\}_{k \le 2})$ for pairwise distinct $X_0, X_1, X_2 \in \XX_n(h) \cap Q_{2\t + 1}(\bar x_{j'})$. In particular, introducing the event
$$E(x_0, x_1, x_2) := \Big\{(r(x_0, x_1), r(x_0, x_1, x_2))\in \Ebd, r(\{x_i\}_{i \le 2}) > r(\{x_i\}_{i \le 1})\Big\},$$
the definition of the factorial moment measures~\eqref{eq:fac_mom} and condition \eqref{mom_eq} yield that
\begin{align*}
	&\P\big(E(X_0, X_1, X_2)
	\text{ for some }X_0, X_1, X_2 \su \XX_n(h) \cap  Q_{2\t + 1}(\bar x_{j '})\big)\\
	&\quad\le c\int_{Q_{2\t + 1}(\bar x_{j '})^2}\one_{\{E(x_0, x_1, x_2)\}}\r_{h, 2, n}(x_0, x_1, x_2)\d(x_0, x_1, x_2)\\
	&\quad\le cC_{\r, 2}\int_{Q_{2\t + 1}(\bar x_{j '})^2} \one_{\{E(x_0, x_1, x_2)\}}\d(x_0, x_1, x_2).
\end{align*}
Noting that by \cite[Proposition 6]{krebs}, the final integral is $O(|\Eb|^{3/4}|\Ed|^{3/4})$ concludes the proof of \eqref{d15_eq}.
\enp

%
%
\bep[Proof of Lemma \ref{covBoundLem}]
     In order to establish the covariance bounds, we proceed along the lines of \cite[Lemma 2]{krebs}. Nevertheless, we present some details since the latter is stated for a cylindrical domain whereas we work with a 3D window growing in two space directions.

First, we introduce independent copies $\{\PP_{i, n}\}_{i \le k_n}$ of $\PP_n$ in order to derive an alternative representation of the martingale differences $\mc D_j$, $j \le k_n$. More precisely, we set
        \begin{align*}
                \PP_{i, n}^* &:= \big(\PP_n \cap \cup_{z \le_{\ms{lex}}z_i}Q(z)\big) \cup \big(\PP_{i, n} \cap \cup_{z >_{\ms{lex}}z_i}Q(z)\big) , \\
                \PP_{i, n}^{**} &:= \big(\PP_n \cap \cup_{z <_{\ms{lex}}z_i}Q(z)\big) \cup \big(\PP_{i, n} \cap \cup_{z \ge_{\ms{lex}}z_i}Q(z)\big) , \\
		X_u^* &:= \prod_{i \in I_u}\big(\b_n(E, \PP_{i, n}^*) - \b_n(E, \PP_{i, n}^{**})\big),
        \end{align*}
	so that $\Cov(X_1,X_2) = \Cov(X_1^*,X_2^*)$.	To ease notation, we henceforth put $R_{i, n}^*:= R(z_i, n; \PP_{i, n}^*)$, $R_{i, n}^{**}:= R(z_i, n; \PP_{i, n}^{**})$ and $R_{i,n}^\vee := R_{i, n}^* \vee R_{i, n}^{**}$. Moreover, in order to separate the influence of the changes at the spatial locations corresponding to the indices in $I_1$ and $I_2$, we also introduce the event
	$$F_n := \Big\{\max_{i \in I_1\cup I_2}R_{i, n}^\vee <\tf14 \dist(\{z_i\}_{i \in I_1},\{z_j\}_{j \in I_2})\Big\}.$$

Next, we express the covariance as
        $
                \Cov(X^*_1, X^*_2) = \Cov\big(X^*_1, X^*_2\one_{\{F_n\}}\big)                + \Cov\big(X^*_1, X^*_2\one_{\{F_n^c\}}\big),
                $
		noting that the random variables $X^*_1$ and $X^*_2\one\{F_n\}$ are independent by condition \ref{S1}. Thus, the first covariance vanishes. For the second one, we deduce from the Cauchy-Schwarz inequality that 
        $\Cov\big(X^*_1, X^*_2\one_{\{F_n^c\}}\big) \le \sqrt{\Var\big(X^*_1\big)}\sqrt{\Var\big(X^*_2\one_{\{F_n^c\}}\big)}.$
        A further application of the Cauchy-Schwarz inequality yields that
	$$\Var\big(X^*_2\one_{\{F_n^c\}}\big)\le \sqrt{\E[(X^*_2)^4}]\sqrt{\P(F_n^c)}.$$ 
	Since the exponential stabilization implies stretched exponential decay of $\P(F_n^c)$, the assertion follows.
\enp

\section{Auxiliary results}
\label{aux_sec}
In this section, we provide the proofs of three auxiliary results appearing in Sections \ref{vor_sec} and \ref{md_sec}.

\begin{figure}[!ht]\centering
  \definecolor{qqqqcc}{rgb}{0.,0.,0.8}
\definecolor{xdxdff}{rgb}{0.49019607843137253,0.49019607843137253,1.}
\definecolor{qqqqff}{rgb}{0.,0.,1.}
\definecolor{zzttqq}{rgb}{0.6,0.2,0.}
\definecolor{cqcqcq}{rgb}{0.7529411764705882,0.7529411764705882,0.7529411764705882}
\begin{tikzpicture}[scale=0.5,line cap=round,line join=round,>=triangle 45,x=1.0cm,y=1.0cm]
\clip(-5.283243740197932,-4.711173003082594) rectangle (10.059409442431445,7.059796657833544);
\fill[color=zzttqq,fill=zzttqq,fill opacity=0.1] (0.,3.) -- (0.,-2.) -- (5.,-2.) -- (5.,3.) -- cycle;
\draw [color=zzttqq] (0.,3.)-- (0.,-2.);
\draw [color=zzttqq] (0.,-2.)-- (5.,-2.);
\draw [color=zzttqq] (5.,-2.)-- (5.,3.);
\draw [color=zzttqq] (5.,3.)-- (0.,3.);
\draw (-4.,6.)-- (1.4365658969228272,1.025293385971551);
\draw (1.4365658969228272,1.025293385971551)-- (4.844061435292844,6.412562868422478);
\draw (1.208130441836569,-0.36435563246985436)-- (-4.521792223243739,-2.5725650316370188);
\draw (1.208130441836569,-0.36435563246985436)-- (3.6257390081661343,-3.4863068519820524);
\draw (1.,-4.)-- (1.,-2.);
\draw (2.,5.)-- (2.,3.);
\draw (4.90236006705965,5.727256503163703) node[anchor=north west] {$C_{i,\,n}$};
\draw (-1.9518933535233347,-1.5065329079011462) node[anchor=north west] {$C_{j,\,n}$};
\draw [color=qqqqcc](1.1317473365421015,-3.4482342761343427) node[anchor=north west] {$G_{j,\,n}(h)$};
\draw [color=qqqqcc](0.3503683954356174,6.156167865448057) node[anchor=north west] {$G_{i,n}(h)$};
\begin{scriptsize}
\draw[color=zzttqq] (2.540670596506408,0.6636039154183084) node {$A$};
\draw [fill=qqqqff] (2.,0.) circle (2.5pt);
\draw[color=qqqqff] (2.4409085501054563,0.13998702071277576) node {$x$};
\draw [fill=black] (-4.,6.) circle (2.5pt);
\draw [fill=black] (1.4365658969228272,1.025293385971551) circle (2.5pt);
\draw [fill=black] (4.844061435292844,6.412562868422478) circle (2.5pt);
\draw [fill=black] (1.208130441836569,-0.36435563246985436) circle (2.5pt);
\draw [fill=black] (-4.521792223243739,-2.5725650316370188) circle (2.5pt);
\draw [fill=black] (3.6257390081661343,-3.4863068519820524) circle (2.5pt);
\draw [fill=qqqqff] (1.,-4.) circle (2.5pt);
\draw [fill=qqqqff] (1.,-2.) circle (2.5pt);
\draw[color=qqqqff] (1.1700578659888594,-1.6017143475204205) node {$Q_j$};
\draw [fill=qqqqff] (2.,5.) circle (2.5pt);
\draw [fill=xdxdff] (2.,3.) circle (2.5pt);
\draw[color=xdxdff] (2.178981125953166,3.404829376453409) node {$Q_i$};
\end{scriptsize}
\end{tikzpicture}
\caption{Sketch for the proof of Lemma~\ref{lem:cells}.}
\label{proof_cells}
 \end{figure}
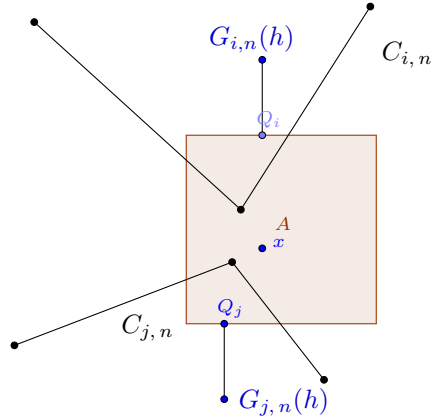

\bel[Voronoi Lemma]
\label{lem:cells}
Let $x \in \R^p$ and $\vp \su \R^p$ be locally finite such that $x + [-1/2, 1/2]^p$ is hit by $\ell \ge 1$ cells of the Voronoi tessellation on $\vp$. Then there exists some $m \ge 2p + 1$ such that ${\vp\big(B(x,\,m) \sm B(x,\,{m - 2p - 1})\big)}  \ge  \ell$.
\enl
\bep
For a visualization of the proof situation see Figure~\ref{proof_cells}. Assume $\{P_1, \dots, P_\ell\} \su \vp$ are points whose Voronoi cells with respect to $\vp$ intersect $x + [-1/2, 1/2]^p$ and let $a := |P_1 - x|$ denote the Euclidean distance from $P_1$ to $x$. Then, the proof of the lemma is completed with $m := \lceil a\rceil + p$ once we show that 
$$\{P_1, \dots, P_\ell\} \su B(x, a + p) \sm B(x, a - p).$$
To prove this claim, fix $1 < j \le \ell$ and let $Q, Q' \in x + [-1/2, 1/2]^p$ denote points that are contained in the Voronoi cells of $P_1$, $P_j$, respectively. In particular,
$$|P_j - x| \ge |P_j - Q| - |Q - x| \ge |P_1 - Q| - |Q - x| \ge a - \sqrt p> m-2p - 1$$
(note that $a - \sqrt p$ can be negative, in which case $B(x,\,m-2p-1)=\emptyset$) and 
$$|P_j - x| \le |P_j - Q'| + |Q' - x| \le |P_1 - Q'| + |Q' - x| \le a + \sqrt p<m.$$

\enp

Recall that 
$$E_{r, x, n} := \bigcap_{\sqrt z \in A_{r, x, n}}\{Q_{\sqrt r}(\sqrt r z) \cap \PP \ne \es\}, \quad \text{where } A_{r, x, n} :=\sqrt r \Z^{d - 1} \cap Q_n \cap Q_r(x)\sm Q_{r - 8 \sqrt r}(x). $$
\bel[$E_{r, x, n}$ occurs with high probability]
\label{lem:er}
There exists $c > 0$ such that for all $r > 1$,
$$\sup_{n \ge 0}\sup_{x \in Q_n} \P(E_{r, x, n}^c) \le \exp(-c r^{(p-1)/2})$$
\enl
\bep
The key idea is to combine the void probabilities for the Poisson point process with a union bound. Indeed, note that $|A_{r, x, n}| \le (\sqrt r + 1)^p$, so that 
$$\P(E_{r, x, n}^c) \le |A_{r, x, n}| \P\big(Q_{\sqrt r}(\sqrt r z) \cap \PP \ne \es\big)  \le C (\sqrt r + 1)^p\exp\big(- r^{(p-1)/2}\big).$$
We conclude the proof since the right-hand side is in $O(\exp(-r^{(p-1)/2}))$.
\enp

We are finally going to show an expression for the cumulant of multiple random variables. Define for a set of indices $S$ the quantities
\[
M^{S}:=\E\Big[\prod_{j\in S}X_j\Big]
\]
and for two subsets $S_1\subset S,\,T_1\subset T$ of indices
\[
U^{S_1,\,T_1}:=\E\Big[\prod_{j\in S_1\cup T_1}X_j\Big]-M^{S_1}M^{T_1}.
\] 
By abuse of notation we write $c(S,\,T)$ in place of $c((X_i)_{i\in S},\,(X_j)_{j\in T})$.
\bel[Proof of~\eqref{eq:fourth_cum}]\label{lemma:fourth_cum} Let $S,\,T$ be a non-trivial partition of $\{1,\,2,\,\dots,\,k\}$. Then,
\[
	c(S,\,T)=\sum_{(S_1,\,T_1),\,(S_2,\,T_2), \dots, (S_r, T_r)}a_{(S_1,\,T_1),\,(S_2,\,T_2), \dots, (S_r, T_r)}U^{S_1,\,T_1}M^{S_2}M^{T_2}\cdots M^{S_r}M^{T_r},
\]
where $a_{(S_1,\,T_1),\,(S_2,\,T_2), \dots, (S_r, T_r)}$ is an integer-valued prefactor, and where the sum runs over all partitions $(S_1, T_1) \cup \cdots \cup (S_r, T_r)$ of $\{1, \dots, k\}$ such that $S_1, T_1$ are non-empty and such that $S_i \su S$, $T_i \su T$.
\enl
\bep
We call ``semi-moment'' of a partition $S_1,\,S_2\subset S$, $T_1,\,T_2\subset T$ the quantity
\[
U^{S_1,\,T_1}M^{S_2}M^{T_2}.
\]
Let $(S,\,T)$ be as in the assumptions. When $k=2$, $S=\{1\}$, $T=\{2\}$ (without loss of generality) we have
\[
c(X_1,\,X_2)=\E[X_1 X_2]-\E[X_1]\E[X_2]=U^{S,\,T}.
\]
Let us show the result for $k>2$. We write~\cite[Eq.~(3.2.6)]{peccatitaqqu}
\[
M^{S\cup T}:=\sum_{(S_1,\,T_1),\,\ldots,\,(S_p,\,T_p)}\prod_{i=1}^p c(S_i,\, T_i)
\]
where $(S_i,\,T_i)_{i=1}^p$ is a partition of $S\cup T$ with each $S_i\subset S,\,T_i\subset T$.  The partitions for which the empty set appears in $(S_i,\,T_i)$ for all $i$ are called ``degenerate'', otherwise they are called ``proper''. Therefore we can split
\begin{equation}\label{eq:split_mom}
M^{S\cup T}=\sum_{(S_1,\,T_1),\,\ldots,\,(S_p,\,T_p)\text{ proper}}\,\prod_{i=1}^p c(S_i,\, T_i)+\sum_{(S_1,\,T_1),\,\ldots,\,(S_p,\,T_p)\text{ degenerate}}\,\prod_{i=1}^p c(S_i,\, T_i).
\end{equation}
We can see that
\[
\sum_{(S_j,\,T_j)_j\text{ degenerate}}\,\prod_{i=1}^p c(S_i,\, T_i)=M^S M^T.
\]
As for the second summand of~\eqref{eq:split_mom}, one can deduce that
\[
\sum_{(S_1,\,T_1),\,\ldots,\,(S_p,\,T_p)\text{ proper}}\,\prod_{i=1}^p c(S_i,\, T_i)=c(S,\, T)+\mathfrak{L} 
\]
where $\mathfrak L$ is a linear combination of semi-moments by the inductive hypothesis. Therefore
\[
c(S,\, T)=M^{S\cup T}-M^S M^T+\mathfrak{L}
\]
where $\mathfrak{L}$ is yet again a (different) linear combination of semi-moments. This concludes the proof of the first statement. To obtain~\eqref{eq:fourth_cum} one can take $k=4$ and $S=\{1\},\,T=\{2,\,3,\,4\}$.
\enp

\end{document}